# THE SOLUTIONS OF THE 3RD AND 4TH CLAY MILLENNIUM PROBLEMS

## The first about the P vs NP in computational complexity and the second about the Navier-Stokes equations

**Konstantinos E. Kyritsis**

## PROLOGUE

In this treatise we present the solutions of the 3rd Clay Millennium problem P vs NP in the Computational Complexity and a proposed solution for the 4th Clay Millennium problem in classical fluid dynamics about the Navier-Stokes equations.

Some initial but incorrect solutions of the 3rd Clay Millennium problem about P vs NP had already been published by me in *International Journal of Pure and Applied Mathematics Volume 120 No. 3 2018, pp 497-510 ISSN: 1311-8080 (printed version); ISSN: 1314-3395 (on-line version) url: http://www.ijpam.eu doi: 10.12732/ijpam.v120i3.1*

But also, in:

Kyritsis C. *On the solution of the 3rd Clay Millennium problem. A short and elegant proof that P ≠ NP in the context of deterministic Turing machines and Zermelo-Frankel set theory.* Proceedings of the first ICQSBEI 2017 conference, Athens, Greece, https://books.google.gr/books?id=BSUsDwAAQBAJ&pg pp 170-181

And also, in:

Kyritsis K. Review of the Solutions of the Clay Millennium Problem about P ≠ NP =EXPTIME World Journal of Research and Review (WJRR) ISSN:2455-3956, Volume-13, Issue-3, September 2021 Pages 21-26

But also in Chapter 3 a 3rd and drastically shorter solution which was presented in the in the 6rth International conference on quantitative, social, biomedical and economic issues, ICQSBE 2022 1st July 2022,

https://icqsbei2022.blogspot.com/2022/06/blog-post.html

http://books.google.com/books/about?id=xZnCEAAAQBAJ

Eventually after discussions and lectures about P vs NP in the School of electrical and Computer engineering in the National Technical University of Athens in 2023, I corrected the solutions to two new solutions presented here

https://www.researchgate.net/publication/397411645_The_solution_of_the_P_vs_NP_millennium_problem_The_non-provability_of_P_is_equal_to_NP_in_1st_order_countable_infinite_logic_of_the_Peano_Arithmetic_based_on_the_Rice-Shapiro_theorem

The proposed solutions of the 4thClay Millennium problem about the Navier-Stokes equations contain a genuine a solution of the initial formulation, based on the CKN theory of blow-ups and an argument by reduction to contradiction. (That there cannot be a velocity blow-up)

Before the final such proposed solution of this millennium problem there was also still another partial solution based again on the hypotheses of conservation of particles here

Kyritsis, K. November 2017 "On the 4th Clay Millennium problem: Proof of the regularity of the solutions of the Euler and Navier-Stokes equations, based on the conservation of particles" *Journal of Scientific Research and Studies Vol 4 (11) , pp304-317,November 2017.*

It seems that at the beginning of each century has become a tradition to state a list of significant and usually difficult problems in the mathematics, that it is considered that their solution will advance significantly the mathematical sciences. At the begging of the 20th century (1900) it was D. Hilbert who formulated and listed 23 problems that most of them have been solved till today (see e.g. https://en.wikipedia.org/wiki/Hilbert%27s_problems) . Those problems from the 23 that have been solved proved to be indeed critical for the overall evolution of mathematics and their applications. Continuing this tradition, the Clay Mathematical Instituted formulated in 2000, 7 critical problems and this time there is a monetary award for their solution (see e.g. http://www.claymath.org/millennium-problems) . From them, the 6th problem (Poincare Hypothesis) it has been accepted that it has been solved by Grigoriy Perelman in 2003. It is not presented here a common or joint method of solution of the 3rd and 4th Clay millenniums problems. It is only because I am an interdisciplinary researcher that I have worked, on both of them. And of course, I had both the advantages and disadvantages of an interdisciplinary researcher. The disadvantage was that I had to sharpen by specialized knowledge in two different areas of Computer science and Mathematical physics , that specialist would not need not do it, while the advantage , that turned out to be more important, were that “I was not blinded by the trees so as to see the forest”; In other words I used new heuristic methods from other disciplines to discover the correct direction of solutions and afterwards I worked out a standard classical proof for each one of them. This is well known in the history of mathematics. E.g. Archimedes found at first the correct formulae of volumes of the sphere, cylinder etc with water, sand and balanced moments of forces experiments before he worked out logically complete proofs of them in the context of Euclidean geometry. Similarly, Newton discovered the laws of gravitation for earth, sun, moon etc. with his, at that time unpublished calculus of fluxes or infinitesimals, and then worked strict proofs within Euclidean geometry in his famous Principia Mathematica.

Similarly, I used myself a heuristic methodology based on statistical mechanics and the particle structure of fluids. Unfortunately, the mathematical models of the fluid dynamic within which this Millennium problem has been formulated are based on the concept of **infinite divisible matter** (before the discovery in the science of physics that matter consists from finite many atoms) and this is a main source of difficulty involving this problem.

Both problems had at least two different directions of solution. For the 3rd Clay Millennium problem, it is:

1) that the non-deterministic polynomial complexity symbolized by NP is equal to a polynomial complexity symbolized by P (in which case the usual setting of passwords and messages might be unsafe) or
2) to a higher e.g. EXPTIME (in which case the usual setting of passwords and messages is as expected to be safe). The heuristic analysis gave that it should hold NP=EXPTIME, which was eventually proved in two different ways.

And for the 4th Clay Millennium problem two different directions of solution would be that:

1) There exist a Blow-up of velocities in finite time.
2) No blow-up exist in finite time and the solutions of the Navier-Stokes equations are regular.

The heuristic analysis gave that because of finite initial energy and energy conservation and mainly because of incompressibility there cannot be a velocity

Blow-up which was eventually proved within the context of classical fluid dynamics that allows for infinite limits etc. More on the logic and strategy of proof for each problem in the next two parts of this treatise.

## TABLE OF CONTENTS

### PART A.

### PART B

**CHAPTER 6**

PROPOSED SOLUTION OF THE MILLENNIUM PROBLEM ON THE NAVIER-STOKES EQUATIONS BASED ON THE PHYSICAL PRINCIPLE OF PARTICLES STRUCTURE CONSERVATION OF THE FLUIDS.

.......................................................................................................................................

**CHAPTER 7**

UNDER THE HYPOTHESES OF THE INITIAL CONDITIONS OF THE MILLENNIUM PROBLEM ABOUT THE NAVIER-STOKES EQUATIONS, THERE CANNOT BE SHAPED IN FINITE TIME, ANY VELOCITIES BLOW-UPS TO INFINITE
(SINGULARITIES) .......................................................................................................
.............................

# PART A.

## CHAPTER 1

## THE SOLUTION OF THE MILLENNIUM PROBLEM ABOUT THE P vs NP IN COMPUTATIONAL COMPLEXITY. INFORMAL DISCUSSION

### Prologue

The standard formulation of the 3rd Clay Millennium problem can be found in (**Cook, Stephen** ***April 2000*** *The P versus NP Problem (PDF), Clay Mathematics Institute site*. http://www.claymath.org/millennium-problems/p-vs-np-problem

http://www.claymath.org/sites/default/files/pvsnp.pdf*)*

1) The P versus NP is a difficult problem, that has troubled the scientific community for some decades

2) It may have simple proofs of a few paragraphs, hopefully not longer than the proof of the Time Hierarchy theorem, which seems to be a deeper result.

3) But it can also have very lengthily and complex proofs, that may take dozens of pages.

What the final proof in the next published is or is not:

1) It does not introduce new theoretical concepts in computational complexity theory so as to solve the P versus NP.

2) It does not use relativization and oracles

3) It does not use diagonalization arguments, although the main proof, utilizes results from the time hierarchy theorem

4) It is not based on improvements of previous bounds of complexity on circuits

5) It is proved with the method of counter-example. Thus, it is transparent short and "simple". It takes any Exptime-complete DTM decision problem, and from it, it derives in the context of deterministic Turing machines a decision problem language which it is apparent that it belongs in the NP class decision problems while it does not belong the class P of decision problems.

6) It seems a "simple" proof because it chooses the right context to make the arguments and constructions and the key-abstraction mentioned above. So, it helps that the scientific community will accept that this 3rd Clay Millennium problem has already been solved.

## RS0) The disruptive role for the computational complexity of the solution of the millennium problem P vs NP. The Russell's disallowed impredicative prdicates and contradictions.

## My odyssey when I tried to re-solve the mystery of P vs NP in the context of computational complexity.

(In the symbolism of paragraphs we utilize RS meaning resolution)

## 0. INTRODUCTION. MY INITIAL APPROACH IN MY ODYSSEY TO SOLVE THE P VS NP PROBLEM.

The P vs NP problem (polynomial time versus non-deterministic polynomial time) si one of the major problems of computational complexity , and the 3rd millennium problem (of the Clay mathematical Institute https://www.claymath.org/). I will not spend space here in formulating it and explain about it. The non-expert reader is refered to the excellent presentation by S Cook in [18] . The P versus NP problem is generally considered not yet solved and by the more careful expert researchers as not yet been known if it has been solved or not. Many have claimed solutions from 2000 to 2016. G. J. Woeginger (see G. J. Woeginger [1] and [2] Wikipedia) compiled a list of more than 100 purported proofs of P = NP from 1986 to 2016, of which 50 were proofs of P ≠ NP, 2 were proofs the problem is unprovable, and one was a proof that it is undecidable I myself read some very few short solutions from this list which turned out to be incorrect. But I did not analyse the long ones (some more than 70 pages). I do not know any researcher who has gone through all the solutions in the list of G. J. Woeginger, to find which, if any, of the more than 100 solutions is correct. Although for some solutions in this list, it is easy to prove that are not correct, no-one has ever published any proof that all of them are not correct. Probably this should be the task of the Clay Mathematical Institute which sponsored the formulation of this problem as one of the 7 Millennium problems, in other words to hire a group of experts to do this task. Nevertheless, according to the rules about the millennium problems of the Clay Mathematical Institute, the Institute is waiting for the community of relevant experts and researchers to indicate by citations if there is a correct solution to the P vs. NP problem. Most of these more than 100 solutions are not in the main journals of complexity theory and the reason is that the most widely read such journals avoid refereeing any solution of the P vs. NP problem for obvious or non-obvious reasons, except perhaps if it is from a very well-known and celebrated professor in the field of complexity. Therefore, there is an obvious social barrier to publishing solutions to this problem in relevant Journals that are widely read. Strangely enough, the monetary award for the solution to this problem had two opposite effects. First an increased number of researchers from all areas trying to solve it, and second an avoidance of the main Journals in the specialization area, to consider solutions to this problem for refereeing which of course would exclude correct solutions also.

In the history of mathematics, it is known that difficult problems that have troubled a lot the mathematicians turned out to have different proofs one simple and one very complex. Such an example is if the general 5th-order polynomial equation can be solved with addition, subtraction, multiplication, division and extraction of radicals starting from the coefficients. The famous mathematician Niels Henrik Abel gave a very simple proof, of not more than 5 pages. On the other hand, the proof of the same, by the E. Galois theory, is a whole book of dozens of pages!

And a famous mathematician once said that "*Once a proof is known to a mathematical problem*, *then immediately after it becomes trivial*!"

It is important to mention, a statement, that is usually attributed to the famous mathematician Yuri Manin, that "A correct proof in mathematics is considered a proof only if it has passed the social barrier of being accepted and understood by the scientific community and published in accepted Journals".

*Passing the obstruction of the social barrier*, *sometimes is more difficult than solving the mathematical problem itself*!

It is similar to the solution of the P versus NP problem in this paper.

The P vs. NP is not a problem that a computer experiment can decide, but rather a problem that requires the correct arguments over the relevant concepts. It is in theoretical computational complexity which utilizes concepts like, "languages of infinite many words", and the infinite is not existing in the computer practice (on the contrary some computer practitioners may consider it a computer worm!).So when I started studying the P vs. NP problem, the first that I asked myself was, **"from which axioms, should I start reasoning?"** Soon I realized that I should start reasoning from the axioms of the mathematical set theory. But this is not enough either. One must determine the size and type of the formal logic allowed so as to have a possible and correct informal proof.

Therefore the P versus NP problem is in fact a set of different problems when they are in the context of different axiomatic systems of set theory and different types and sizes of logic (e.g. $1^{st}$ order countable logic, higher order countable logic, higher order uncountable logic etc or in $1^{st}$ order Peano Arithmetic versus $1^{st}$ order ZFC set theory) .

When I finished solving the Millennium problem about the Navier-Stokes equations in fluid dynamics, in 2017 (See [6] ) I started trying to solve also the P vs NP millennium problem. I was afraid that it would be more difficult compared to that in fluid dynamics, because fluid dynamics is a centuries old specialization of mathematics and it has proved practically all the necessary tools about it. On the other hand Computational Complexity was not really more than 50-years old, and it was expected that it has not proved yet all its necessary tools and results that are essential for this specialization. But I never expected my Odyssey that followed and the disruptive role that the solution of the millennium problem P vs NP would have for the common perception of Computational Complexity.

Initially in 2017 I thought that I proved that P is not equal to NP (see e.g. [3], [4] , [5] .[7] , [8].)

***I tried 3 times to prove that P is not equal to NP, but on all 3 times the proof was incorrect. This made me change my strategy, and find the correct goal , which was to trye to prove that P=NP or that P not equal to NP are not provable in 1st order Peano Arithmetic. This resulted in succedng to prove that P=NP is not provable in 1st order Arithmetic. Initially I thought I proved also that P not equal to NP is also non-provable in 1st order Peano Arithmetic. But I was wrong. Actually I did manage to prove the P is not equal to NP, but not in 1st order Peano Arithmetic but in ZFC set theory. If it is alsoprovable in 1st order Peano Arithmetic or not is in general a different issue.***

So I spent all my sabbatical of 2023-2024 in trying to re-solve it.

From my experience in solving the Clay millennium problem in fluid dynamics about the Navier-Stokes equations (See references [7] ), I knew that a probable failure in solving it is that the mainstream of research has forgotten, or disregarded as too old research, some old but significant researchers who had accumulated significant results and theorem relevant to the subject. In the case of the Navier-Stokes equations it was the invariants discovered by Helmholtz, Kelvin and Stokes. As I am an interdisciplinary researcher, is a pleasure to me to search all possible research even distantly relevant to this millennium problem. I wanted to be sure, that I did not miss any results that could be used.

## ABOUT THE LEVELS ⊢ , ⊨ , M⊨ , OF LOGICAL TRUTH AND PROOF BY REDUCTION TO CONTADICTION.

The proof under discussion, is of the type of **REDUCTION TO CONTRADICTION** (εις άτοπό απαγωγή, in ancient Greek) so allow me to refresh and clarify, why this method is a weaker proof than the other methods and what we know about it . In symbolic logic of formal systems, there are 3 concepts of "truth" about **propositions.**

1) **Proof from the axioms** (see e.g. references [105] page 163), which is usually denoted by ⊢

2) **Validity in any interpretation** (by models of set theory) usually denoted by ⊨

3) **Truth value equal to 1 in the model M**, usually for a proposition P denoted by T(P)=1 or we may denote it by M ⊨ P , In a fixed Model M of the theory. For any proposition p either p or ~p are true.

Unfortunately, all these difficulties in the mathematics and logic, occur because we allow the existence in the **ontology of the infinite. This is so because axiomatic theories with finite models are easily provable consistent.** That is why my last part of my research, is the creation of the **Democritus-Pythagoras-Archimedes geometry and differential and integral calculus without the infinite** (only finite many (some trillion many though) points and numbers with finite many digits. I will come in to that again when I will define for you, (see in the attached file) this that all the working computer scientists want, which the evolution of the Church thesis to the concept of **REAL MACHINE** versus the **TURING MACHINE**.

The Real machine has **nothing infinite**.

**A REALMACHINE**

1) It is a Turing machine,

2) with finite bounded writing strip (space bound Smax, for all machines),

3) finite bounded number of states (code bounded Cmax for all machines),

4) finite bounded time to run Tmax (max number-complexity for all machines)

5) and finite bounded number of input words as input data languages (Dmax , for all machines).

**Still historically the concept of the infinite in the mathematical ontology has served thinkers to keep a distance from the material reality and thus being able to think freely.**

Since we mentioned that the assignment of 2-valued truth values to the propositions of an axiomatic system is not unique, and that the proof by reduction to contradiction, proves only the truth value, then we may speculate, that there **in the same axiomatic system may exist different proofs by reduction to contradiction, which prove different assignments of 2-valued truth to propositions. Still we will accept only one as valid, because we assume that we are in a single universal assignment of 2-valued truth, although we don't know which one (actually we chose it partly and gradually as we accumulate theorems) !**

Having discussed the above, we should mention that the famous **P vs NP problem, is not really a single problem**, because we may be interested e.g. to establish **5 factors that will differentiate the problem.**

**3 Logical levels of truth factors**

1) Axiomatic system $\vdash$ P equal or not equal to NP

2) Axiomatic system $\models$ P equal or not equal to NP

3) N $\models$ P equal or not equal to NP (in a particular model e.g. the stanrd model N of the natural numbers )

4) We must specify also , in which **axiomatic system** we are. Are we in the **Neumann-Bernays-Goedel NBG set theory** which allows classes, or are we in the **Zermelo-Frankel ZFC set theory** which does not allow classes, but only sets. Or are we in the $1^{st}$ order logic axiomatic system of the **Peano Arithmetic PA** ?

5) Furthermore beyond all the above, because in the definition of NP a logical symbols **R(x,y)** appears of a binary relation of words (a word x and a certificate y of it) appears, **we must know the Logic in which R is a predicate**. Is a **countable $1^{st}$ order logic** of a set theory? Or is it the countable $1^{st}$ order logic of Peano Arithmetic? Not all countable logics would give the same result, because the predicates of such a logic being countable if at all will define countably many only sets, while the set theory has much more than countable. Or is it an **uncountable informal logic** where all possible sets can appear as predicates, e.g. $2^{nd}$ order logic of Peano Arithmetic in full semantics? Here we restrict to countable $1^{st}$ order logic of PA or of ZFC.

## CHAPTER 2-3

# The solution of the P vs NP. Why the P=NP cannot be proved in 1st order Peano Arithmetic and why eventually, we can prove that P≠NP in ZFC set theory.

**By Kyritsis Konstantinos**
University of Ioannina , School of Economics. ckiritsi@uoi.gr

**Abstract**

The millennium problem P vs NP problem, has resisted any solution more than half a century now. This is not without a good reason. In this work we prove that P ≠NP in ZFC set theory. In order to prove it, we use the method of David Martin in the solution of the 10th Hilbert problem about the Diophantine equations, where he derives a non-provability in the logic of Peano axiomatic system (which is that there is a Diophantine equation which has no solution, but this is not provable), from a non-recursive enumerability of non-solvable Diophantine equations. We transfer this method in the case of the non-recursive enumerability, that the Rice-Shapiro theorem or other methods gives. Then we derive easily the above mentioned result.

## 1. Introduction

The P vs NP problem of computational complexity, is considered as one of the greatest problem of Computational Complexity, in the last half of a century. It is asked if it is true that if a problem that a non-deterministic Turing machines decides in polynomial time, can be also decided by a deterministic Turing machine in polynomial time as well (or not). We shall not spend space here describing the detailed formulation of the problem. The reader is referred to [6] by S. Cook, where the problem is formulated in detail. This problem is also known as the 3rd Millennium problem.

**The difference of**

**a) PA , or ZFC ⊦ P equal or not equal to NP from**

**b) N ⊨ P equal or not equal to NP”**

In this article, we must involve aspects of symbolic logic, and models of formal axiomatic systems. Inparticular we msut disciminate between

1) PA or ZFC ⊢ P equal or not equal to NP

2) N ⊨ P equal or not equal to NP (Truth inside the standard model of Peano Arithmetic PA inside the ZFC set theory.

We assume that all the theory of computation and complexity is carried out inside that standard model N of Peano Arithmetic PA. The standard model N with its truth is defined in the ZFC or NBG set theory.

**2.The simple proof that P≠ NP in 1st order ZFC set theory. .**

To prove this we shall be based on a technique by **Martin Davis** in the Theorem 7.7. in [8] page 263 in his account of the solution of the 10th Hilbert problem where he **is producing non-provability in 1st order countable logic starting with non-recursive enumerability and undecidability of some sets of the objective ontology.**

To allow the reader to have a more clear and direct perception of the Martin Davis method , we copy paste his theorem 7.7 at his exposition [8] of the solution of the 10nth Hilbert problem, together with the short proof by Martin Davis.

**THEOREM 7.7. (MARTIN DAVIS 1973)** *Corresponding to any given axiomatization of number theory, there is a Diophantine equation which has no positive integer solutions, but such that this fact cannot be proved within the given axiomatization.*

**Proof:** A rigorous proof would involve a precise definition of “axiomatization of number theory” which is outside the scope of this article. An informal heuristic argument follows:

One uses the given axiomatization to systematically generate all of the theorems (i.e., consequences of the axioms). Among these theorems will be some asserting that some Diophantine equation has no solution. Whenever such is encountered it is placed on a special list called LIST A. At the same time a list, LIST B, is made of Diophantine equations which have solutions. LIST B is constructed by a search procedure, e.g., at the nth stage of the search look at the first n Diophantine equations (in a suitable list) and test for solutions in which each argument is ≤ n. Thus every Diophantine equation which has positive integer solutions will eventually be placed in LIST B. If likewise each Diophantine equation with no solutions would eventually appear in LIST A, then one would have an algorithm for Hilbert’s tenth problem. Namely, to test a given equation for possession of a solution simply begin generating LIST A and LIST B until the given equation appears in one list or the other. Since Hilbert’s tenth problem is unsolvable, some equation with no solution must be omitted from LIST A. But this is just the assertion of the theorem. **QED.**

**Remark 2.0** The reader may find Martin Davis proof, too intuitive or very informal and it is so. But any one can ask any advanced cloud AI, to enrich this proof with formal details, and make it a normal explicit, and formal proof. *Especially there is the point that in order to separate the axioms of Davis List A, from all the recursive enumerable axioms of PA, he must prove them a decidable set*. But this not difficult, as the syntactical structure of theorems of the type "**The Diophantine equation D, has no solution"** ***is a characteristic repetitive string among all possible theorems.*** We will use now this idea by Martin Davis, but instead of Diophantine equations D (as string of parameters and variables) we have Turing machine M, and their word-codes **#<M>** , where **#<>** is e.g. a **Goedel arithmetization function.**

Let us give some definitions

**Definitions 2.1** *Let any non-trivial language L. We define as set*
***M(L)****={ #<M>/ M is a Turing machine, that accepts the language L}*
*And* ***L(M)*** *is the language accepted by the Turing machine M, with code #<M>.*
*Let us define, some sets of Turing machines from the classical complexity classes:*

***P(M)*** *= { #<M> : L(M) ∈ P } = { <M> : ∃ k ∀x M halts on x in* $\leq |x|^k$ *steps and accepts iff x∈L(M) }*

***NP(M)*** *= { #<M> : ∃ poly-time NTM Turring Machine . N with L(N)=L(M) } = { #<M> : L(M) ∈ NP }*

It is easily realized that M(L), P(M), NP(M) and if P ≠NM , NP(M)-P(M)≠∅ , are non-trivial and semantic properties of Turing machines , thus from the Rice theorem, they are undecidable.

But we want something more, we want to prove that they are is also **not recursive enumerable.** So we have the next lemma

**Lemma 2.0** *The sets of Turing machines codes* M(L), P(M), NP(M) and if P ≠NM then NP(M)-P(M)≠ ∅ *as codes,* ***are not recursive enumerable.***

**Proof:** We can use either the Trakhtenbrot theorem or the Arithmetic Hierarchy or even simpler The **Rice-Shapiro theorem (**see **[24] )** . We do not need to spend space for this here. The reader may ask any advanced cloud AI to present him with the details of such a proof.

**QED.**

**Remark 2.1** The Rice-Shapiro theorem is equivalent with the Turing theorem (undecidability of halting) , which in its turn is equivalent with the $1^{st}$ incompleteness theorem for PA. The details of such an enlightening proof can be given by any advanced cloud AI.

Then we proceed to apply the above mentioned method by Martin Davis, and for this we need to define a sentences **OPT** (for optimistic case and not the pessimistic case) that would be equivalent to the sentences P=NP at least meta-mathematically. So

**OPT corresponds to P=NP**

And here is the sentence

**Definition 2.2** *OPT="For any Turing machine M , with its code #<M> if #<M> belongs to* NP(M) *, the L(M) also belongs to P"*

The standard literature gives us that the sentence OPT, is a $\Pi_2^0$ sentence of the arithmetic hierarchy of formulae. We may therefore put it as an additional axiom of the 1st order Peano Arithmetic axiomatic system and get the PA*=PA+OPT.

**Lemma 2.1** *If the sentence OPT is put as an axiom in the PA axiomatic system, to derive a new one PA*=PA+OPT. If PA* is consistent then by the logical rule of* ***universal instantiation*** *, also the individual sentences*

*For OPT: " If the code #<M> belongs to NP(M) , then L(M) belongs to P" are also correspondingly to OPT , provable one by one.*

**Proof:** Again direct, elementary by the logical rule of **universal instantiation** .
**QED.**

We now restrict from any possible code #<M>, to only the codes of the set *NP(M) ,* for which we know by Lemma 2.0 that the y are not r.e. And similarly we restrict to the set of individual sentences
*" the code #<M> belongs to NP(M) , and L(M) belongs to P"* We prove the next lemma . In the case of the code #<M> belonging to NP(M) such a sentence is equivalent to the sentence *" if the code #<M> belongs to NP(M) , then L(M) belongs to P"*

**Lemma 2.2** *At least one of the sentences " the code #<M> belongs to NP(M) , and L(M) belongs to P" is not a provable theorem.*

**Proof:** Now if each of these sentences were a provable theorem, we would obviously get from its characteristic syntactic string of repetitive pattern, a decidable subset from the r.e. set of all theorems Th(PA) of Peano Arithmetic. Exactly as in the case of Theorem 7.7 by Martin Davis.
Being a decidable subset of a r.e. set Th(PA) it would make the corresponding codes #<M> of NP(M) r.e. contradiction by lemma 2.0. Thus at least one of the above sentences is not a provable theorem **QED.**

We may notice that the recursive enumerability of the relevant theorems is at the meta-mathematical level, **external** to the Peano arithmetic, compared to the non-recursive enumerability of the codes of the Turing machines as by lemma 2.0 which **internal** . It is the same as in the Theorem 7.7 by Martin Davis. But we may make them **both internal in ZFC set theory** by using an appropriate Godel arithmetization which projects the Peano arithmetic internally in ZFC. Then both concets of recursive enumerability are internal in an assumed consistent axiomatic system, the ZFC.

But then this result of Lemma 2.2 contradicts the result of lemma 2.1, thus the premise hypothesis in the lemma 2.1 of the consistency of PA*=PA+OPT must not hold. So we have the next two corollary

**Corollary 2.0.** *If we accept the 1st incompleteness for PA, as a true meta-mathematical claim then: The sentence OPT (corresponding to P=NP) makes the axiomatic system PA*=PA+OPT inconsistent in 1st order logic, therefore it cannot be proved in 1st order Peano Arithmetic. The proof is meta-mathematical for PA and inside the ZFC set theory*

**Proof:** We just gave it before the statement of the corollary, The mentioning of the 1st incompleteness follows from the remark 2.1. **QED**

**Corollary 2.1.** *If we accept the 1st incompleteness for PA, as a true meta-mathematical claim then:The sentence OPT (corresponding to P=NP) makes the axiomatic system PA*=PA+OPT inconsistent in 1st order logic. The proof is meta-mathematical for PA and inside the ZFC set theory. We may conclude that P≠NP is proved and true in N in ZFC set theory.*

**Proof:** The mentioning of the 1st incompleteness follows from the remark 2.1. Because the proof leading to a contradiction was carried out using r.e. of theorems of PA, Th(PA) , and Th(PA*), that are assumed projected inside ZFC set theory (e.g. on to ordinals larger than ω, outside N) we realize that it was a proof outside PA and some parts of it including the contradiction may non-expressible in PA, but definetely expressible inside ZFC. Therefore it is much natural, to consider the argument as a proof by reduction to contradiction (RAA=reductio ad absurdum) inside ZFC rather than inside PA. This makes safer to claim that we proved P≠NP in ZFC and not in PA. We already know that Presburger arithmetic PA(+) is strictly weaker that Peano arithmetic PA=PA(+,*) which in its turn is strictly weaker that PA(+,*,#)=PA(#) , where # is a Goedel arithmetization, which essentially in ZFC. Theoretically someone could prove P≠NP in ZFC, while this might not be possible in PA. An example of a similar situation is the arithmetic sentence is the Goodstein theorem (See [27]) which has been proved as non-provable in PA, while it has been proved in ZFC. QED.

**Remark 2.2.** The author has met with many contradictions arising from the use of vicious circle self-referential diagonal definitions of type II, in the proofs of some of the important theorems of Computational Complexity. (The reader is advised to consult the [22] [23] [24] [25], [26]) These self-referential definitions are proved to make the ZFC set theory contradictory (in [22]). But even without valid proofs we could take these important theorems as axioms of Computational Complexity.

**References.**

## PART B

## CHAPTER 4

## THE PROPOSED SOLUTONS OF THE MILLENNIUM PROBLEM ABOUT THE NAVIER STOKES EQUATIONS

**Prologue.**

The standard formulation of the 4th Clay Millennium problem can be found in the site of the Clay Mathematical Institute here:

http://www.claymath.org/millenniumproblems/navier%E2%80%93stokes-equation and here http://www.claymath.org/sites/default/files/navierstokes.pdf

Roughly speaking it asks if in classical 3 dimensional incompressible fluids , (governed by the Navier-Stokes equations) with finite initial energy and smooth initial conditions (with pressures and velocities falling to zero faster than all polynomial powers as we go to infinite distances away or in short smooth Schwartz initial conditions) the flow will continuous forever smooth or would there be a finite time, where velocities and pressures will blow-up to infinite and smoothness will break? The standard formulation is both with periodic initial conditions or not periodic.

Most of the mathematicians were expecting that, since it has been proved that there is no blow-up in 2-dimensions, this should hold in 3 dimensions too. But as more than half a century has passed with not being able to prove it many researchers started believing that because of the vortex stretching which is possible only in 3-dimasions and not in 2-dimensions a blow-up might exist.

Because it was easier to do at the beginning, I spent about half a year discovering more than a dozen of explicitly formulated cases of axial symmetric flows that lead to blow-up in finite time. Nevertheless, for all of them, it was necessary that they start with infinite initial energy and the initial vorticities were unbounded.

So I went back to the more probable case that no Blow-up can occur in finite time.

My heuristic analysis which took 1-2 years, with statistical mechanics and classical fluid dynamics in digital differential and integral calculus suggested to me that there should not exist in finite time a blow-up. The naïve and simple argument was that a blow up would give that at least one particle of the fluid (and in statistical mechanics or classical fluid dynamics in digital differential and integral calculus, finite many finite particles do exist) would exhibit infinite kinetic energy. Nevertheless, what is easy to prove in heuristic context is not at all easy to prove in the classical context of fluid dynamics where there are not finite many particles of finite and lower bounded size, but infinite many points with zero size.

In this strategy my interdisciplinary approach was an advantage. I did not consider as consistent for sciences that e.g. statistical mechanics would give that there is no-blow up in finite time, while classical fluid dynamics would prove that there is a blow-up in finite time.

The next table makes the comparisons in statistical mechanics and classical fluid dynamics

**Table 0**

| COMPARISON AND MUTUAL SIGNIFICANCE OF DIFFERENT TYPES OF MATHEMATICAL MODELS FOR THE 4TH CLAY PROBLEM (NO EXTERNAL FORCE) | **CONTINUOUS FLUID MECHANICS MODEL** | **STATISTICAL MECHANICS MODEL** |
|---|---|---|
| SMOOTH SCHWARTZ INITIAL CONDITIONS | YES | POSSIBLE TO IMPOSE |
| FINITE INITIAL ENERGY | YES | YES |
| CONSERVATION OF THE PARTICES | YES(NON-OBVIOUS FORMULATION) | YES (OBVIOUS FORMULATION) |
| LOCAL SMOOTH EVOLUTION IN A INITIAL FINITE TIME INTERVAL | YES | POSSIBLE TO DERIVE |
| EMERGENCE OF A BLOW-UP IN FINITE TIME | IMPOSSIBLE TO OCCUR | IMPOSSIBLE TO OCCUR |

So as it was easy to prove in statistical mechanics that there is no blow-up in finite time, I thought , so as to increase our confidence for the correct side of the solution of the problem , to add hypotheses to the standard formulation of the 4$^{th}$ Clay Millennium problem that correspond to the conservation of particles during the flow, and which would lead to an accessible solution of this problem (that there is no Blow-up in finite) dew to finite initial energy and energy and particle conservation. This of course was not the solution of the 4$^{th}$ Clay Millennium problem, and the solution finally is presented in the 2$^{nd}$ and last paper in this part B of this treatise in chapter 6.

in the chapter 6 the millennium problem is solved again under a simpler hypothesis about the particle structure conservation of the fluid which is based on the vorticity, and so the proof is even simpler.

In the next I have created **3 different solutions** of the millennium problem about the Navier-Stokes equations. The first two , Chapter 5, 6, are based on formulating the physical fact of existence of particles and particles conservation in the fluid, which is an extra hypothesis, to the original formulation, thus not complete solution of the millennium problem. But the 3rd solution on chapter 7, is based only on the CKN theory of blow-ups and does not introduce any additional hypothesis to the original formulation of the problem. Thus the 3rd solution is a genuine solution of the 4th Millennium problem. In this 3rd solution we utilize the existing CKN theory of possible blow-ups, in particular, type I and type II, and we prove by reduction to contradiction that neither of the two types of blow-ups can appear in finite time blow-ups because this would mean that a blow-up or singularity will happen in a 3D region too. The proof is based on the non-compressible of the fluid, which implies that the 1D integrals on paths of the pressure forces is path-independent (conservative vector field) and that the initial energy is finite.

**What is that maybe we do not understand very well with the Navier-Stokes equations and their physics?**

In statistical mechanical models of incompressible flow, we have the realistic advantage of finite many particles, e.g. like balls B(x,r) with finite diameter r. These particles as they flow in time, remain particles of the same nature and size and the velocities and inside them remain approximately constant.

E. g. for the case water, we are speaking here for molecules of $H_2O$, that are estimated to have a diameter of 2.75 angstroms or 2r= 2.75*10^(-10) meters. The Representative Elementary Volume or REV is usually assumed about **$10^{-15}$ meters.** But if we go deaper e.g. for electrons, then the size grows smaller between **$10^{-23}$ meters** and $10^{-15}$ meters. Furthermore all the electrons , protons and neutrons are assumed to spin with about the same frequency which is about 1 terahertz. Or as **vorticity** ω it would be 2π times that or **ω=6.26 Terraherz.** This frequency is what makes our molecular matter physical reality coherent, and it make increase or decrease simultaneously in large areas from deeper factors, that might as well be related to human consciousness too, as some experiments prove.

Fluid dynamics from the strictly mathematical point of view goes smaller than any scale, as it was conceived at the centuries that matter was believed to be infinite divisible. It is somehow a shame that we still use such type of models! In practice when we apply e.g. Navier-Stokes fluid dynamics in water fluids we assume that we may go as small as the Representative Elementary Volume or REV is usually assumed about **$10^{-15}$ meters.** In the principle of the particle structure conservation that we define and use in this article, we assume that the vorticity ω of the **Representative Elementary Volumes but the point vorticities too are uniformly modulated and influences by the constant spin of ω=6.26 Terraherz of the electrons, neutrons and protons.**

Because space and time dimensions in classical fluid dynamics goes in orders of smallness , smaller and at least as small as the real physical molecules , atoms and particles of the fluids, this might suggest imposing too, such conditions resembling uniform continuity conditions. In the case of continuous fluid dynamics models such natural conditions, emerging from the particle nature of material fluids, together with the energy conservation, the incompressibility and the momentum conservation, as laws conserved in time, may derive the regularity of the local smooth solutions of the Navier-Stokes equations.

In [8] it is defined **a concept of particle structure and conservation of particle structure** in incomprehensible fluids, which is used to propose a solution of the Millennium problem about the Navier-Stokes eqautions in favor of the regularity (no blow-up). But here we define a simple concept of particle structure and conservation of particle structure, based only on the **vorticity ,** which gives as simpler solution of the millennium problem about the Navier-Stokes equations in favor of the regularity.

**Definition 3.0 The physical principle of particle structure conservation**

*Let a non-compressble fluid, which satisfies the Euler or Navier-Stokes equations, and flows smoothly on all of 3D space and in the time interval [0, T). it may be defined as non-periodic or periodic flow, and without forces. We say that it has at time t=0* ***particle structure*** *relative to vorticity* ***at smallness scale ε*** *iff for any point x, where it has vorticity ω(x,0), there is a diffeomorphic image B(0) of a spherical ball, with volume equal to |B(0)|=ε , which contains x , and the 3D volume average of the vorticity on B(0) is equal to ω(x,0). In symbols*

$$\int_{B(0)} \omega dv = \omega(x,0)$$

*We say that we have* ***volume preserving particle structure conservation*** *during [0,T) , if the same happens with the* ***same ε,*** *for all later times t in [0,T).*

Obviously here we may take as ε the volume of a Representative Elementary Volume or REV is usually assumed about **$10^{-15}$ meters,** thus ($10^{-15}$ meters)$^3$=ε or

**ε=10^(-45) m$^{3}$.**

Such a physical principle may seem very strong and unusual, to the standards of fluid dynamics, that pretends that the matter is infiniote divisible, but it has physical anchors. To see that **it is compatible with the Navier-Stokes equations** , we only need to notice that the irrotational flows (or potential flows) do satisfy it, and obviously from the physical point of view, it must be applied for a smallness scale ε, much larger than the ε=10^(-45) $m^{3}$.

It is also important to realize that , we define volume preserving particle structure conservation, thus **we must have non-compressible volume preserving flows. For compressible fluids this definition does not apply.**

Besides the forgotten conservation law of finite particles, there are **two more forgotten laws of conservation or invariants**. The first of them is the obvious that during the flow, the physical measuring units dimensions (dimensional analysis) of the involved physical quantities (mass density, velocity, vorticity, momentum, energy, force point density, pressure, etc.) are conserved. It is not very wise to eliminate the physical magnitudes interpretation and their dimensional analysis when trying to solve the millennium problem, because the dimensional analysis is a very simple and powerful interlink of the involved quantities and leads with the physical interpretation, to a transcendental shortcut to symbolic calculations. **By eliminating the dimensional analysis, we lose part of the map to reach our goal.**

The 3rd forgotten conservation law or invariant, is related to the viscosity (friction). Because we do know that at each point (pointwise), the viscosity is only subtracting kinetic energy, with an irreversible way, and converting it to thermal energy, (negative energy point density), and this is preserved in the flow, (it can never convert thermal energy to macroscopic kinetic energy), we know that its sign does not change too, it is a flow invariant, so the integrated 1D or 2D work density is always of the same sign (negative) and as sign, an invariant of the flow. The **conservation or invariance of the sign of work density by the viscosity (friction)** is summarized in the lemma 3.1 below.

**Viscosity irreversibility Principle 3.2 The absolute physical principle of irreversibility of the dissipated kinetic energy by viscosity, for non-compressible fluids as forgotten invariant.**

*For incompressible viscous fluids, the dissipated energy by viscosity, only subtracts from the kinetic energy and convert it to thermal. It is not reversible, it cannot add to the kinetic energy. And this applies not only to the total 3D space and time integral, but also to the 2D and 1D densities of rate of energy dissipation (absolute irreversibility=both global as 3D integral but also local as 2D, 1D or point values densities of rates) . As a consequence, on a trajectory path, although the*

*viscosity force density is not aligned with the velocity, its projection on the trajectory path, always opposes the pressure force densities projection on the trajectory path. Similarly, at the circulation of velocities multiplied by the constant density (as a 1D circular momentum density), the work density of the viscosity force densities on the same circle, will always subtract kinetic energy and reduce the momentum density. If the mathematics of solutions of the NSE equations, allow violations of the above irreversibility of the dissipation energy, then they are not acceptable as physical meaningful solutions. In the same way that solutions where the kinetic energy is increasing are not acceptable as physically meaningful solutions (as L. Fefferman remarks in [3] in the formulation of the millennium problem).*

Although we shall use the physical principle of volume preserving particle structure conservation, we shall not use in the next in this article the physical principle of irreversibility of the dissipated kinetic energy.

## CHAPTER 5

## ON THE SOLUTION OF THE 4TH MILLENNIUM PROBLEM. PROOF OF THE REGULARITY OF THE SOLUTIONS OF THE EULER AND NAVIER-STOKES EQUATIONS, BASED ON THE CONSERVATION OF PARTICLES.

### Abstract

As more and more researchers tend to believe that with the hypotheses of the official formulation of the 4th Clay Millennium problem a blowup may occur, a new goal is set: to find the simplest and most physically natural enhancement of the hypotheses in the official formulation so that the regularity can be proved in the case of 3 dimensions too. The position of this paper is that the standard assumptions of the official formulation of the 4th Clay millennium problem, although they reflect, the finiteness and the conservation of momentum and energy and the smoothness of the incompressible physical flows, they do not reflect the conservation of particles as local structure. By formulating the later conservation and adding it to the hypotheses, we prove the regularity (global in time existence and smoothness) both for the Euler and the Navier-Stokes equations.

### 1. Introduction

This 1st paper is an initial version of the published paper in the Journal of Scientific Research and Studies Vol. 4(11), pp. 304-317, November, 2017 ISSN 2375-8791

The famous problem of the 4th Clay mathematical Institute as formulated in FEFFERMAN C. L. 2006 , is considered a significant challenge to the science of mathematical physics of fluids, not only because it has withstand the efforts of the scientific community for decades to prove it (or types of converses to it) but also because it is supposed to hide a significant missing perception about the nature of our mathematical formulations of the physical flows through the Euler and the Navier-Stokes equations.

When the 4th Clay Millennium problem was officially formulated the majority was hoping that the regularity was holding also in 3 dimensions as it had been proved to hold also in 2 dimensions. But as time passed more and more mathematicians started believing that a Blowup can occur with the hypotheses of the official formulation. Therefore, a new goal is set to find the simplest and most physically natural enhancement of the hypotheses in the official formulation so that the regularity can be proved in the case of 3 dimensions too. This is done by the current paper.

After 3 years of research, in the 4th Clay Millennium problem, the author came to believe that, what most of the mathematicians would want, (and seemingly including the official formulators of the problem too), in other words a proof of the regularity in 3 dimensions as well, cannot be given merely by the assumptions of the official formulation of the problem. In other words, a Blow-up may occur even with compact support smooth initial data with finite energy. But solving the 4th Clay Millennium problem, by designing such a case of Blow-up is I think not interesting from the physical point of view, as it is quite away from physical applications and a mathematical pathological curiosity. On the other hand, discovering what physical aspect of the flows is not captured by the mathematical hypotheses, is I believe a more significant contribution to the science of mathematical physics in this area. Although the mathematical assumptions of the official formulation reflect, the finiteness and the conservation of momentum and energy and the smoothness of the incompressible physical flows, they do not reflect the conservation of particles as local structure. By adding this physical aspect formulated simply in the context of continuous fluid mechanics, the expected result of regularity can be proved.

In statistical mechanical models of incompressible flow, we have the realistic advantage of finite many particles, e.g. like balls B(x,r) with finite diameter r. These particles as they flow in time, **remain particles of the same nature and size** and the velocities and inside them remain approximately constant.

Because space and time dimensions in classical fluid dynamics goes in orders of smallness, smaller and at least as small as the real physical molecules, atoms and particles of the fluids, this might suggest imposing too, such conditions resembling uniform continuity conditions. In the case of continuous fluid dynamics models such natural conditions, emerging from the particle nature of material fluids, together with the energy conservation, the incompressibility and the momentum conservation, as laws conserved in time, may derive the regularity of the local smooth solutions of the Euler and Navier-Stokes equations. **For every atom or material particle of a material fluid, we may assume around it a ball of fixed radius, called _particle range_ depending on the size of the atom or particle, that covers the particle and** **a little bit of the electromagnetic, gravitational or quantum vacuum field around it, that their velocities and space-time accelerations are affected by the motion of the molecule or particle.** E.g. for the case water, we are speaking here for molecules of $H_2O$, that are estimated to have a diameter of 2.75 angstroms or 2r= 2.75*10^(-10) meters, we may define as water molecule **particle range** the balls $B(r_0)$ of radius $r_0$=3*10^(-10) meters around the water molecule. As the fluid flows, especially in our case here of incompressible fluids, the shape and size of the molecules do not change much, neither there are significant differences of the velocities and space-time accelerations of parts of the molecule. Bounds $\delta_u$ $\delta_\omega$ of such differences remain constant as the fluid flows. **We may call this effect as the *<u>principle of conservation of particles</u>* as a local structure.** This principle must be posed in equal setting as the energy conservation and incompressibility together with the Navier-Stokes or Euler equations. Of course, if the fluid is say of solar plasma matter, such a description would not apply. Nevertheless, then incompressibility is hardly a property of it. But if we are talking about incompressible fluids that the molecule is conserved as well as the atoms and do not change atomic number (as e.g. in fusion or fission) then this principle is physically valid. The principle of conservation of particles as a local structure, blocks the self-similarity effects of concentrating the energy and turbulence in very small areas and creating thus a Blow-up. It is the missing invariant in the discussion of many researchers about supercritical, critical and subcritical invariants in scale transformations of the solutions.

The exact definition of the conservation of particles as local structure Is in DEFINITION 5.1 and it is as follows:

(**<u>Conservation of particles as local structure in a fluid</u>**)

*Let a smooth solution of the Euler or Navier-Stokes equations for incompressible fluids, that exists in the time interval [0,T). We may assume initial data on all of $R^3$ or only on a connected compact support $V_0$. For simplicity let us concentrate only on the latter simpler case. Let us denote by F the displacement transformation of the flow. Let us also denote by g the partial derivatives of 1st order in space and time, that is $\left|\partial_x^a \partial_t^b u\ (x)\right|$, $|\alpha|=1$, $|b|<=1$, and call then space-time accelerations . We say that there is **<u>conservation of the particles</u>** **<u>in the interval [0,T) in</u>** a derivatives homogenous setting, as a local structure of the solution if and only if:*

*There is a small radius r, and small constants $\delta_x$ , $\delta_u$ , $\delta_\omega$ , >0 so that for all t in [0,T) there is a finite cover $C_t$ (in the case of initial data on $R^3$ · it is infinite cover, but finite on any compact subset) of $V_t$ , from balls B(r) of radius r, called **<u>ranges of the particles</u>** , such that:*

*1) For an $x_1$ and $x_2$ in a ball B(r) of $V_s$ , s in [0,T), ||F($x_1$ )-F($x_2$)||<=r+ $\delta_x$ for all t>=s in [0,T).*

*2) For an $x_1$ and $x_2$ in a ball B(r) of $V_s$ ,s in [0,T), ||u(F($x_1$ ))-u(F($x_2$))||<= $\delta_u$ for all t >=s in [0,T).*

*3) For an $x_1$ and $x_2$ in a ball B(r) of $V_s$ , s in [0,T), ||g(F($x_1$ ))-g(F($x_2$))||<= $\delta_\omega$ for all t >=s in [0,T).*

*If we state the same conditions 1) 2) 3) for all times t in [0,+∞) , then we say that we have the **strong version** of the conservation of particles as local structure.*

We prove in paragraph 5 in PROPOSITION 5.2 that indeed adding the above conservation of particles as local structure in the hypotheses of the official formulation of the 4th Clay Millennium problem, we solve it, in the sense of proving the regularity (global in time smoothness) of the locally in time smooth solutions that are known to exist.

A short outline of the logical structure of the paper is the next.

1) The paragraph 3, contains the official formulation of the 4th Cay millennium problem as in FEFFERMAN C. L. 2006. The official formulation is any one of 4 different conjectures, that two of them, assert the existence of blow-up in the periodic and non-periodic case, and two of them the non-existence of blow-up , that is the global in time regularity in the periodic and non-periodic case. We concentrate on to prove the regularity in the non-periodic case or conjecture (A) with is described by equations 1-6 after adding the conservation of particles as a local structure. The paragraph 3 contains definitions, and more modern symbolism introduced by T, Tao in TAO T. 2013. The current paper follows the formal and mathematical austerity standards that the official formulation has set, together with the suggested by the official formulation relevant results in the literature like in the book MAJDA A.J-BERTOZZI A. L. 2002.

But we try also not to lose the intuition of the physical interpretation, as we are in the area of mathematical physics rather than pure mathematics.

The goal is that reader after reading a dozen of mathematical propositions and their proofs, he must be able at the end to have simple physical intuition, why the conjecture (A) of the 4th Clay millennium together with the conservation of particles in the hypotheses problem holds.

2) The paragraph 4 contains some known theorems and results, that are to be used in this paper, so that the reader is not searching them in the literature and can have a direct, at a glance, image of what holds and what is proved. The most important are a list of necessary and sufficient conditions of regularity (PROPOSITIONS 4.5-4.10)

The same paragraph contains also some well-known and very relevant results that are not used directly but are there for a better understanding of the physics.

3) The paragraph 5 contains the main idea that the conservation of particles during the flow can be approximately formulated in the context of continuous fluid mechanics and that is the key missing concept of conservation that acts as subcritical invariant in other words blocks the self-similar concentration of energy and turbulence that would create a Blowup. With this new invariant we prove the regularity in the case of 3 dimensions: PROPOSITIONS 5.2.

**4)** The paragraph 6 contains the idea of defining **a measure of turbulence** in the context of deterministic mechanics based on the **total variation** of the component functions or norms (DEFINITION 6.1) It is also made the significant observation that the smoothness of the solutions of the Euler and Navier-Stokes equations is not a general type of smoothness but one that would deserve the name **homogeneous smoothness** (Remark 6.2) .

According to CONSTANTIN P. 2007 "…The blowup problem for the Euler equations is a major open problem of PDE, theory of far greater physical importance that the blow-up problem of the Navier-Stokes equation, which is of course known to non-specialists because of the Clay Millennium problem…"

Almost all of our proved propositions and in particular the regularity in paragraphs 4 , 5 and 6 (in particular PROPOSITION 4.11 and PROPOSITION 5.2) are stated not only for the Navier-Stokes but also for the Euler equations.

## 2. The ontology of the continuous fluid mechanics models versus the ontology of statistical mechanics models. The main physical idea of the proof of the regularity in 3 spatial dimensions.

All researchers discriminate between the physical reality with its natural physical ontology (e.g. atoms, fluids etc) from the mathematical ontology (e.g. sets, numbers, vector fields etc). If we do not do that much confusion will arise. The main difference of the physical reality ontology, from the mathematical reality ontology, is what the mathematician D. Hilbert had remarked in his writings about the infinite. He remarked that nowhere in the physical reality there is anything infinite, while the mathematical infinite, as formulated in a special axiom of the infinite in G. Cantor's theory of sets, is simply a convenient phenomenological abstraction, at a time that the atomic theory of matter was not well established yet in the mathematical community. In the physical reality ontology, as best captured by statistical mechanics models, the problem of the global 3-dimensional regularity seems easier to solve. For example it is known (See PROPOSITION 4.9 and PROPOSITION 4.12 maximum Cauchy development, and it is referred also in the official formulation of the Clay millennium problem in C. L. FEFFERMAN 2006) that if the global 3D regularity does not hold then the velocities become unbounded or tend in absolute value to infinite as time gets close to the finite Blow-up time. Now we know that a fluid consists from a finite number of atoms and molecules, which also have finite mass and with a lower bound in their size. **If such a phenomenon (Blowup) would occur, it would mean that for at least one particle the kinetic energy, is increasing in an unbounded way.** But from the assumptions (see paragraph 3) the initial energy is finite, so this could never happen. We conclude that the fluid is 3D globally in time regular. Unfortunately, such an argument although valid in statistical mechanics models (see also MURIEL A 2000), in not valid in continuous fluid mechanics models, where there are not atoms or particles with lower bound of finite mass, but only points with zero dimension, and only mass density. We must notice also here that this argument is not likely to be successful if the fluid is compressible. In fact, it has been proved that a blow-up may occur even with smooth compact support initial data, in the case of compressible fluids. One of the reasons is that if there is not lower bound in the density of the fluid, then even without violating the momentum and energy conservation, a density converging to zero may lead to velocities of some points converging to infinite. Nevertheless, if we formulate in the context of continuous fluid mechanics the conservation of particles as a local structure (DEFINITION 5.1) then we can derive a similar argument (see proof of PROPOSITION 5.1) where **if a Blowup occurs in finite time, then the kinetic energy of a finite small ball (called in** DEFINITION 5.1 **particle-range) will become unbounded**, which is again impossible, due to the hypotheses if finite initial energy and energy conservation.

The next table compares the hypotheses and conclusions both in continuous fluid mechanics models and statistical mechanics models of the 4$^{th}$ Clay millennium problem in its officially formulation together with the hypothesis of conservation of particles. It would be paradoxical that we would be able to prove the regularity in statistical mechanics and we would not be able to prove it in continuous fluid mechanics.

**Table 1**

| COMPARISON AND MUTUAL SIGNIFICANCE OF DIFFERENT TYPES OF MATHEMATICAL MODELS FOR THE 4TH CLAY PROBLEM (NO EXTERNAL FORCE) | **CONTINUOUS FLUID MECHANICS MODEL** | **STATISTICAL MECHANICS MODEL** |
|---|---|---|
| SMOOTH SCHWARTZ INITIAL CONDITIONS | YES | POSSIBLE TO IMPOSE |
| FINITE INITIAL ENERGY | YES | YES |
| ***CONSERVATION OF THE PARTICES*** | ***YES(NON-OBVIOUS FORMULATION)*** | ***YES (OBVIOUS FORMULATION)*** |
| LOCAL SMOOTH EVOLUTION IN A INITIAL FINITE TIME INTERVAL | YES | POSSIBLE TO DERIVE |
| EMERGENCE OF A BLOW-UP IN FINITE TIME | IMPOSSIBLE TO OCCUR | IMPOSSIBLE TO OCCUR |

**3. The official formulation of the Clay Mathematical Institute 4$^{th}$ Clay millennium conjecture of 3D regularity and some definitions.**

In this paragraph we highlight the basic parts of the official formulation of the 4$^{th}$ Clay millennium problem, together with some more modern, since 2006, symbolism, by relevant researchers, like T. Tao.

***In this paper I consider the conjecture (A) of C. L. FEFFERMAN 2006 official formulation of the 4$^{th}$ Clay millennium problem, which I identify throughout the paper as the 4$^{th}$ Clay millennium problem.***

The Navier-Stokes equations are given by (by R we denote the field of the real numbers, ν>0 is the viscosity coefficient)

$$\frac{\partial}{\partial t}u_i + \sum_{j=1}^{n} u_j \frac{\partial u_i}{\partial x_j} = -\frac{\partial p}{\partial x_i} + \nu \Delta u_i \qquad (x \in R^3 , t>=0 , n=3) \qquad \text{(eq.1 )}$$

$$divu = \sum_{i=1}^{n} \frac{\partial u_i}{\partial x_i} = 0 \qquad (x \in R^3 , t>=0 , n=3) \qquad \text{(eq.2 )}$$

with initial conditions u(x,0)=$u^0$(x) **xεR³** and $u_0$ (x) C∞ divergence-free vector field on **R³** (eq.3)

$$\Delta = \sum_{i=1}^{n} \frac{\partial^2}{\partial x_i^2}$$ is the Laplacian operator. The Euler equations are when ν=0

For physically meaningful solutions we want to make sure that $u^0$(x) does not grow large as |x|→∞. This is set by defining $u^0$(x) and called in this paper **Schwartz initial conditions**, in other words

$\left|\partial_x^a u^0(x)\right| \le C_{a,K}(1+|x|)^{-K}$ on **R³** for any α and K (eq.4)

(Schwartz used such functions to define the space of Schwartz distributions)

We accept as physical meaningful solutions only if it satisfies

p, u ∈$C^\infty(R^3 \times [0,\infty))$ (eq.5)

and

$\int_{\Re^3} |u(x,t)|^2 dx < C$ for all t>=0 (**Bounded or finite energy**) (eq.6)

The conjecture (A) of he Clay Millennium problem (case of no external force, but homogeneous and regular velocities) claims that for the Navier-Stokes equations, v>0, n=3 , with divergence free , Schwartz initial velocities , there are for all times t>0 , smooth velocity field and pressure, that are solutions of the Navier-Stokes equations with bounded energy**, in other words satisfying the equations eq.1 , eq.2 , eq. 3, eq.4 , eq.5 eq.6** . It is stated in the same formal formulation of the Clay millennium problem by C. L. Fefferman see C. L. FEFFERMAN 2006 (see page 2nd line 5 from below) that the conjecture (A) has been proved to holds locally. "..if the time internal [0,∞), is replaced by a small time interval [0,T), with T depending on the initial data....". In other words there is ∞>T>0, such that there is continuous and smooth solution u(x,t) ∈ $C^\infty(R^3 \times [0,T))$. In this paper, as it is standard almost everywhere, the term smooth refers to the space $C^\infty$

Following TAO, T 2013, we define some specific terminology, about the hypotheses of the Clay millennium problem, that will be used in the next.

*We must notice that the definitions below can apply also to the case of inviscid flows, satisfying the* ***Euler*** *equations.*

DEFINITION 3.1 (Smooth solutions to the Navier-Stokes system). *A smooth set of data* for the Navier-Stokes system up to time T is a triplet ($u_0$, f, T), where 0 < T < ∞ is a time, the initial velocity vector field $u_0 : R^3 \to R^3$ and the forcing term f : [0, T] $\times R^3 \rightarrow R^3$ are assumed to be smooth on $R^3$ and $[0, T] \times R^3$ respectively (thus, $u_0$ is infinitely differentiable in space, and f is infinitely differentiable in space time), and $u_0$ is furthermore required to be divergence-free:

$\nabla \cdot u_0 = 0.$

If f = 0, we say that the data is *homogeneous.*

In the proofs of the main conjecture, we will not consider any external force, thus the data will always be homogeneous. But we will state intermediate propositions with external forcing. Next, we are defining simple diffentiability of the data by Sobolev spaces.

DEFINITION 3.2 We define the $H^1$ norm (or enstrophy norm) $H^1(u_0, f, T)$ of the data to be the quantity

$H^1(u_0, f, T) := \|u_0\|_{H^1_X(R^3)} + \|f\|_{L^\infty_t H^1_X(R^3)} < \infty$ and say that $(u_0, f, T)$ *is $H^1$* if

$H^1(u_0, f, T) < \infty.$

DEFINITION 3.3 We say that a *smooth set of data* $(u_0, f, T)$ is *Schwartz* if, for all integers α, m, k ≥ 0, one has

$$\sup_{x \in R^3} (1+|x|)^k \left|\nabla_x^a u_0(x)\right| < \infty \quad \text{and} \quad \sup_{(t,x)\in[0,T]\times R^3} (1+|x|)^k \left|\nabla_x^a \partial_t^m f(x)\right| < \infty$$

Thus, for instance, the solution or initial data having Schwartz property implies having the $H^1$ property.

DEFINITION 3.4 A *smooth solution* to the Navier-Stokes system, or a *smooth solution* for short, is a quintuplet $(u, p, u_0, f, T)$, where $(u_0, f, T)$ is a *smooth set of data*, and the velocity vector field $u : [0, T] \times R^3 \rightarrow R^3$ and pressure field $p : [0, T]\times R^3 \rightarrow R$ are smooth functions on $[0, T]\times R^3$ that obey the Navier-Stokes equation (eq. 1) but with external forcing term f,

$$\frac{\partial}{\partial t}u_i + \sum_{j=1}^{n} u_j \frac{\partial u_i}{\partial x_j} = -\frac{\partial p}{\partial x_i} + v\Delta u_i + f_i \quad (x \in R^3,\ t>=0,\ n=3)$$

and also the incompressibility property (eq.2) on all of $[0, T] \times R^3$, but also the initial condition $u(0, x) = u_0(x)$ for all $x \in R^3$

DEFINITION 3.5 Similarly, we say that $(u, p, u_0, f, T)$ *is $H^1$* if the associated data $(u_0, f, T)$ *is $H^1$*, and in addition one has

$$\|u\|_{L^\infty_t H^1_X([0,T]\times R^3)} + \|u\|_{L^2_t H^2_X([0,T]\times R^3)} < \infty$$

We say that the solution is *incomplete in [0,T),* if it is defined only in [0,t] for every t<T.

We use here the notation of *mixed norms* (as e.g. in TAO, T 2013). That is if $\|u\|_{H_x^k(\Omega)}$ is the classical Sobolev norm ,of smooth function of a spatial domain Ω, $u:\Omega \to R$ , I is a time interval and $\|u\|_{L_t^p(I)}$ is the classical $L^p$ -norm, then the mixed norm is defined by

$$\|u\|_{L_t^p H_x^k(I\times\Omega)} := (\int_I \|u(t)\|_{H_x^k(\Omega)}^p dt)^{1/p} \text{ and } \|u\|_{L_t^\infty H_x^k(I\times\Omega)} := ess\sup_{t\in I}\|u(t)\|_{H_x^k(\Omega)}$$

Similar instead of the Sobolev norm for other norms of function spaces.

We also denote by $C_x^k(\Omega)$ , for any natural number $k \geq 0$ , the space of all k times continuously differentiable functions $u:\Omega \to R$ , with finite the next norm

$$\|u\|_{C_x^k(\Omega)} := \sum_{j=0}^{k} \|\nabla^j u\|_{L_x^\infty(\Omega)}$$

We use also the next notation for *hybrid norms*. Given two normed spaces X, Y on the same domain (in either space or time), we endow their intersection $X \cap Y$ with the norm

$$\|u\|_{X\cap Y} := \|u\|_X + \|u\|_Y .$$

In particular in the we will use the next notation for intersection functions spaces, and their hybrid norms.

$$X^k(I\times\Omega) := L_t^\infty H_x^k(I\times\Omega) \cap L_x^2 H_x^{k+1}(I\times\Omega) .$$

We also use the *big O notation*, in the standard way, that is X=O(Y) means

$X \leq CY$ for some constant C. If the constant C depends on a parameter s, we denote it by $C_s$ and we write $X=O_s(Y)$.

We denote the difference of two sets A, B by A\B. And we denote Euclidean balls by $B(a,r) := \{x \in R^3 : |x-a| \leq r\}$ , where |x| is the Euclidean norm.

With the above terminology the target Clay millennium conjecture in this paper can be restated as the next proposition

**The 4th Clay millennium problem (Conjecture A)**

(**Global regularity for homogeneous Schwartz data**). *Let (*$u_0$*, 0, T) be a homogeneous Schwartz set of data. Then there exists a smooth finite energy solution (u, p,* $u_0$*, 0, T) with the indicated data (notice it is for any T>0, thus global in time)*

.

## 4. Some known or directly derivable, useful results that will be used.

In this paragraph I state, some known theorems and results, that are to be used in this paper, so that the reader is not searching them in the literature and can have a direct, at a glance, image of what holds and what is proved.

A review of this paragraph is as follows:

Propositions 4.1, 4.2 are mainly about the uniqueness and existence locally of smooth solutions of the Navier-Stokes and Euler equations with smooth Schwartz initial data. Proposition 4.3 are necessary or sufficient or necessary and sufficient conditions of regularity (global in time smoothness) for the Euler equations without viscosity. Equations 8-15 are forms of the energy conservation and finiteness of the energy loss in viscosity or energy dissipation. Equations 16-18 relate quantities for the conditions of regularity. Proposition 4.4 is the equivalence of smooth Schwartz initial data with smooth compact support initial data for the formulation of the 4th Clay millennium problem. Propositions 4.5-4.9 are necessary and sufficient conditions for regularity, either for the Euler or Navier-Stokes equations, while Propositions 4.10 is a necessary and sufficient condition of regularity for only the Navier-Stokes with non-zero viscidity.

In the next I want to use, the basic local existence and uniqueness of smooth solutions to the Navier-Stokes (and Euler) equations, that is usually referred also as the well posedness, as it corresponds to the existence and uniqueness of the physical reality causality of the flow. The theory of well-posedness for smooth solutions is summarized in an adequate form for this paper by the Theorem 5.4 in TAO, T. 2013.

I give first the definition of **mild solution** as in TAO, T. 2013 page 9. Mild solutions must satisfy a condition on the pressure given by the velocities. Solutions of smooth initial Schwartz data are always mild, but the concept of mild solutions is a generalization to apply for non-fast decaying in space initial data , as the Schwartz data, but for which data we may want also to have local existence and uniqueness of solutions.

**DEFINITION 4.1**

We define a *$H^1$ mild solution* (u, p, $u_0$, f, T) to be fields u, f :[0, T] × $R^3$ → $R^3$,

p : :[0, T] × $R^3$ → R, $u_0$ : $R^3$ → $R^3$, with 0 < T < ∞ , obeying the regularity hypotheses

$$u_0 \in H^1_x(R^3)$$

$$f \in L^\infty_t H^1_x([0,T]\times R^3)$$

$$u \in L^\infty_t H^1_x \cap L^2_t H^2_x([0,T]\times R^3)$$

with the pressure p being given by (Poisson)

$$p = -\Delta^{-1}\partial_i\partial_j(u_i u_j) + \Delta^{-1}\nabla\cdot f \qquad \text{(eq. 7)}$$

(Here the summation conventions is used , to not write the Greek big Sigma).

which obey the incompressibility conditions (eq. 2), (eq. 3) and satisfy the integral form of the Navier-Stokes equations

$$u(t) = e^{t\Delta}u_0 + \int_0^t e^{(t-t')\Delta}(-(u\cdot\nabla)u - \nabla p + f)(t')dt'$$

with initial conditions u(x,0)=$u^0$(x) **.**

We notice that the definition holds also for the in viscid flows, satisfying the Euler equations. The viscosity coefficient here has been normalized to ν=1**.**

In reviewing the local well-posedness theory of $H^1$ mild solutions, the next can be said. The content of the theorem 5.4 in TAO, T. 2013 (that I also state here for the convenience of the reader and from which derive our PROPOSITION 4.2) is largely standard (and in many cases it has been improved by more powerful current well-posedness theory). I mention here for example the relevant research by PRODI G 1959 and SERRIN, J 1963, The local existence theory follows from the work of KATO, T. PONCE, G. 1988, the regularity of mild solutions follows from the work of LADYZHENSKAYA, O. A. 1967. There are now a number of advanced local well-posedness results at regularity, especially that of KOCH, H., TATARU, D.2001.

There are many other papers and authors that have proved the local existence and uniqueness of smooth solutions with different methods. As it is referred in C. L. FEFFERMAN 2006 I refer to the reader to the MAJDA A.J-BERTOZZI A. L. 2002 page 104 Theorem 3.4,

I state here for the convenience of the reader the summarizing theorem 5.4 as in TAO T. 2013. I omit the part (v) of Lipchitz stability of the solutions from the statement of the theorem. I use the standard O() notation here, x=O(y) meaning x<=cy for some absolute constant c. If the constant c depends on a parameter k, we set it as index of $O_k()$.

It is important to remark here that the existence and uniqueness results locally in time (well-posedness) , hold also not only for the case of viscous flows following the Navier-Stokes equations, but also for the case of inviscid flows under the Euler equations. There are many other papers and authors that have proved the local existence and uniqueness of smooth solutions both for the Navier-Stokes and the Euler equation with the same methodology, where the value of the viscosity coefficient v=0, can as well be included. I refer e.g. the reader to the MAJDA A.J-BERTOZZI A. L. 2002-page 104 Theorem 3.4, paragraph 3.2.3, and paragraph 4.1 page 138.

**PROPOSITION 4.1** *(Local well-posedness in $H^1$). Let ($u_0$, f, T) be $H^1$ data.*

*(i) (Strong solution) If (u, p, $u_0$, f, T) is an $H^1$ mild solution, then*

$$u \in C_t^0 H_x^1([0,T]\times R^3)$$

*(ii) (Local existence and regularity) If*

$$(\|u_0\|_{H_X^1(R^3)} + \|f\|_{L_t^1 H_X^1(R^3)})^4 T < c$$

*for a sufficiently small absolute constant c > 0, then there exists*

*a $H^1$ mild solution (u, p, $u_0$, f, T) with the indicated data, with*

$$\|u\|_{X^k([0,T]\times R^3)} = O(\|u_0\|_{H_X^1(R^3)} + \|f\|_{L_t^1 H_X^1(R^3)})$$

*and more generally*

$$\|u\|_{X^k([0,T]\times R^3)} = O_k(\|u_0\|_{H_X^k(R^3)}, \|f\|_{L_t^1 H_X^1(R^3)}, 1)$$

*for each k>=1 . In particular, one has local existence whenever*

*T is sufficiently small, depending on the norm $H^1$($u_0$, f, T).*

*(iii) (Uniqueness) There is at most one $H^1$ mild solution (u, p, $u_0$, f, T)*

*with the indicated data.*

*(iv**) (Regularity)** If (u, p, $u_0$, f, T ) is a $H^1$ mild solution, and ($u_0$, f, T)*

*is (smooth) Schwartz data, then u and p is smooth solution; in fact, one has*

$\partial_t^j u, \partial_t^j p \in L_t^\infty H^k([0,T]\times R^3)$ *for all j, K >=0.*

For the proof of the above theorem, the reader is referred to the TAO, T. 2013 theorem 5.4, but also to the papers and books, of the above mentioned other authors**.**

Next I state the local existence and uniqueness of smooth solutions of the Navier-Stokes (and Euler) equations with smooth Schwartz initial conditions, that I will use in this paper, explicitly as a PROPOSITION 4.2 here.

**PROPOSITION 4.2 Local existence and uniqueness of smooth solutions or smooth well posedness**. *Let* $\mathbf{u_0(x)}$ , $\mathbf{p_0(x)}$ *be smooth and Schwartz initial data at t=0 of the Navier-Stokes (or Euler) equations, then there is a finite time interval [0,T] (in general depending on the above initial conditions) so that there is a unique smooth local in time solution of the Navier-Stokes (or Euler) equations*

$\mathbf{u(x)}$ , $\mathbf{p(x)}$ $\in C^\infty(R^3 \times [0,T])$

**Proof**: We simply apply the PROPOSITION 4.1 above and in particular, from the part (ii) and the assumption in the PROPOSITION 4.2, that the initial data are smooth Schwartz , we get the local existence of $H^1$ mild solution (u, p, $u_0$, 0, T). From the part (iv) we get that it is also a smooth solution. From the part (iii), we get that it is unique.

As an alternative we may apply the theorems in MAJDA A.J-BERTOZZI A. L. 2002-page 104 Theorem 3.4, paragraph 3.2.3, and paragraph 4.1 page 138, and getthe local in time solution, then derive from the part (iv) of the PROPOSITION 4.1 above, that they are also in the classical sense smooth. QED.

**Remark 4.1** We remark here that the property of smooth Schwartz initial data, is not in general conserved in later times than t=0, of the smooth solution in the Navier-Stokes equations, because it is a very strong fast decaying property at spatially infinity. But for lower rank derivatives of the velocities (and vorticity) we have the **(global and) local energy estimate**, and **(global and) local enstrophy estimate** theorems that reduce the decaying of the solutions at later times than t=0, at spatially infinite to the decaying of the initial data at spatially infinite. See e.g. TAO, T. 2013, Theorem 8.2 (Remark 8.7) and Theorem 10.1 (Remark 10.6).

Furthermore, in the same paper of formal formulation of the Clay millennium conjecture, L. FEFFERMAN 2006 (see page 3rd line 6 from above), it is stated that the 3D global regularity of such smooth solutions is controlled by the **bounded accumulation in finite time intervals** of the vorticity (Beale-Kato-Majda). I state this also explicitly for the convenience of the reader, for smooth solutions of the Navier-Stokes equations with smooth Schwartz initial conditions, as the **PROPOSITION 4.6 When we say here bounded accumulation** e.g. of the deformations D, **on finite internals**, we mean in the sense e.g. of the proposition 5.1 page 171 in the book MAJDA A.J-BERTOZZI A. L. 2002, which is a definition designed to control the existence or not of finite blowup times. In other words, for any finite time interval

[0, T], there is a constant M such that

$$\int_0^t |D|_{L^\infty}(s)ds <= M$$

I state here for the convenience of the reader, a well-known proposition of equivalent necessary and sufficient conditions of existence globally in time of solutions of the Euler equations, as inviscid smooth flows. It is the proposition 5.1 in MAJDA A.J-BERTOZZI A. L. 2002 page 171.

The *stretching* is defined by

$S(x,t) =: D\xi \cdot \xi$ if $\xi \neq 0$ and $S(x,t) =: 0$ if $\xi = 0$ where $\xi =: \frac{\omega}{|\omega|}$, ω being the vortcity.

**PROPOSITION 4.3** Equivalent *Physical Conditions for Potential Singular Solutions of the Euler equations. The following conditions are equivalent for smooth Schwartz initial data:*

*(1) The time interval,* $[0, T^*)$ *with* $T^* < \infty$ *is a maximal interval of smooth* $H^s$

*existence of solutions for the 3D Euler equations.*

*(2) The vorticity ω accumulates so rapidly in time that*

$$\int_0^t |\omega|_{L^\infty}(s)ds \to +\infty \text{ as } t \text{ tends to } T^*$$

*(3) The deformation matrix D accumulates so rapidly in time that*

$$\int_0^t |D|_{L^\infty}(s)ds \to +\infty \text{ as } t \text{ tends to } T^*$$

*(4) The stretching factor S(**x**, t) accumulates so rapidly in time that*

$$\int_0^t [\max_{x\in R^3} S(x,s)]ds \to +\infty \text{ as } t \text{ tends to } T^*$$

The next theorem establishes the equivalence of smooth connected compact support initial data with the smooth Schwartz initial data, for the homogeneous version of the 4th Clay Millennium problem. It can be stated either for local in time smooth solutions or global in time smooth solutions. The advantage assuming connected compact support smooth initial data, is obvious, as this is preserved in time by smooth functions and also integrations are easier when done on compact connected sets.

**Remark 4.2 Finite initial energy and energy conservation equations:**

When we want to prove that the smoothness in the local in time solutions of the Euler or Navier-Stokes equations is conserved, and that they can be extended indefinitely in time, we usually apply a "reduction ad absurdum" argument: Let the maximum finite time T* and interval [0,T*) so that the local solution can be extended smooth in it.. Then the time T* will be a blow-up time, and if we manage to extend smoothly the solutions on [0,T*]. Then there is no finite Blow-up time T* and the solutions holds in [0,+∞). Below are listed necessary and sufficient conditions for this extension to be possible. Obviously not smoothness assumption can be made for the time T*, as this is what must be proved. But we still can assume that at T* the energy conservation and momentum conservation will hold even for a singularity at T*, as these are universal laws of nature, and the integrals that calculate them, do not require smooth functions but only integrable functions, that may have points of discontinuity.

A very well known form of the energy conservation equation and accumulative energy dissipation is the next:

$$\frac{1}{2}\int_{R^3}\|u(x,T)\|^2 dx+\int_0^T\int_{R^3}\|\nabla u(x,t)\|^2 dxdt=\frac{1}{2}\int_{R^3}\|u(x,0\|^2 dx \qquad \text{(eq. 8 )}$$

Where:

$$E(0)=\frac{1}{2}\int_{R^3}\|u(x,0)\|^2 dx \qquad \text{(eq. 9)}$$

is the initial finite energy

$$E(T)=\frac{1}{2}\int_{R^3}\|u(x,T)\|^2 dx \qquad \text{(eq. 10)}$$

is the final finite energy

and $\Delta E=\int_0^T\int_{R^3}\|\nabla u(x,t)\|^2 dxdt$ (eq. 11)

is the accumulative finite energy dissipation from time 0 to time T , because of viscosity in to internal heat of the fluid. For the Euler equations it is zero. Obviously

ΔE<=E(0)>=E(T) (eq. 12)

The rate of energy dissipation is given by

$$\frac{dE}{dt}(t) = -v\int_{R^3} \|\nabla u\|^2 dx < 0 \qquad \text{(eq. 13)}$$

(v, is the viscosity coefficient. See e.g. MAJDA, A.J-BERTOZZI, A. L. 2002 Proposition 1.13, equation (1.80) pp. 28)

At this point we may discuss, that for the smooth local in time solutions of the Euler equations, in other words for flows without viscosity, it is paradoxical from the physical point of view to assume, that the total accumulative in time energy dissipation is zero while the time or space-point density of the energy dissipation (the former is the $\|\nabla u(x,t)\|^2{}_{L_\infty}$), is not zero! It is indeed from the physical meaningful point of view unnatural, as we cannot assume that there is a loss of energy from to viscosity at a point and a gain from "anti-viscosity" at another point making the total zero. Neither to assume that the time and point density of energy dissipation is non-zero or even infinite at a space point, at a time, or in general at a set of time and space points of measure zero and zero at all other points, which would still make the total accumulative energy dissipation zero. **The reason is of course that the absence of viscosity, occurs at every point and every time, and not only in an accumulative energy level.** If a physical researcher does not accept such inviscid solutions of the Euler equation as having physical meaning, then for all other solutions that have physical meaning and the $\|\nabla u(x,t)\|^2{}_{L_\infty}$ is zero (and come so from appropriate initial data), we may apply the PROPOSITION 4.7 below and **deduce directly, that the local in time smooth solutions of the Euler equations, with smooth Schwartz initial data, and finite initial energy, and zero time and space point energy dissipation density due to viscosity, are also regular (global in time smooth).** For such regular inviscid solutions, we may see from the inequality in (eq. 15) below, that the total L2-norm of the vorticity is not increasing by time. We capture this remark in PROPOSITION 4.11 below.

**Remark 4.3** The next are 3 very useful inequalities for the unique local in time [0,T], smooth solutions u of the Euler and Navier-Stokes equations with smooth Schwartz initial data and finite initial energy (they hold for more general conditions on initial data, but we will not use that):

By $||.||_m$ we denote the Sobolev norm of order m. So, if m=0 itis essentially the $L_2$-norm. By $||.||_{L\infty}$ we denote the supremum norm, u is the velocity, ω is the vorticity, and cm, c are constants.

**1)** $\|u(x,T)\|_m \leq \|u(x,0)\|_m \exp(\int_0^T c_m \|\nabla(u(x,t)\|_{L_\infty} dt)$ (eq. 14)

(see e.g. MAJDA, A.J-BERTOZZI, A. L. 2002, proof of Theorem 3.6 pp117, equation (3.79))

**2)** $\|\omega(x,t)\|_0 \leq \|\omega(x,0)\|_0 \exp(c\int_0^T \|\nabla u(x,t)\|_{L_\infty} dt)$ (eq. 15)

(see e.g. MAJDA, A.J-BERTOZZI, A. L. 2002, proof of Theorem 3.6 pp117, equation (3.80))

**3)** $\|\nabla u(x,t)\|_{L_\infty} \leq \|\nabla u(x,0)\|_0 \exp(\int_0^t \|\omega(x,s)\|_{L_\infty} ds)$ (eq. 16)

(see e.g. MAJDA, A.J-BERTOZZI, A. L. 2002, proof of Theorem 3.6 pp118, last equation of the proof)

The next are a list of well know necessary and sufficient conditions, for regularity (global in time existence and smoothness) of the solutions of Euler and Navier-Stokes equations, under the standard assumption in the $4^{th}$ Clay Millennium problem of smooth Schwartz initial data, that after theorem Proposition 4.4 above can be formulated equivalently with smooth compact connected support data. We denote by T* be the maximum Blow-up time (if it exists) that the local solution u(x,t) is smooth in [0,T*).

**PROPOSITION 4.5 (Necessary and sufficient condition for regularity)**
*The local solution u(x,t) , t in [0,T*) of the Euler or Navier-Stokes equations, with smooth Schwartz initial data, can be extended to [0,T*], where T* is the maximal time that the local solution u(x,t) is smooth in [0,T*), if and only if the* ***Sobolev norm*** *||**u(x,t)**||*$_m$*, m>=3/2+2 , remains bounded , by the same bound in all of [0,T*), then , there is no maximal Blow-up time T*, and the solution exists smooth in [0,+∞)*

**Remark 4.4** See e.g. . MAJDA, A.J-BERTOZZI, A. L. 2002, pp 115, line 10 from below)

**PROPOSITION 4.6 (Necessary and sufficient condition for regularity.** Beale-Kato-Majda**)**

*The local solution u(x,t) , t in [0,T*) of the Euler or Navier-Stokes equations, with smooth compact connected support initial data, can be extended to [0,T*], where T* is the maximal time that the local solution u(x,t) is smooth in [0,T*), if and only if for the finite time interval [0,T*], there exist a bound M>0, so that the* ***vorticity has bounded by M, accumulation*** *in [0,T*]:*

$$\int_0^{T^*} \|\omega(x,t)\|_{L_\infty} dt \leq M \qquad \text{(eq17)}$$

*Then there is no maximal Blow-up time T*, and the solution exists smooth in [0,+∞)*

**Remark 4.5** See e.g. . MAJDA, A.J-BERTOZZI, A. L. 2002, pp 115, Theorem 3.6. Also page 171 theorem 5.1 for the case of inviscid flows. . See also LEMARIE-RIEUSSET P.G. 2002. Conversely if regularity holds, then in any interval from the smoothness in a compact connected set, the vorticity is supremum bounded. The above theorems in the book MAJDA A.J-BERTOZZI A. L. 2002 guarantee that the above conditions extent the local in time solution to global in time, that is to solutions (u, p, $u_0$, f, T) which is $H^1$ mild solution, **for any T**. Then applying the part (iv) of the PROPOSITION 4.1 above, we get that this solution is also smooth in the classical sense, for all T>0, thus globally in time smooth.

**PROPOSITION 4.7 (Necessary and sufficient condition for regularity)**

*The local solution u(x,t) , t in [0,T*) of the Euler or Navier-Stokes equations, with smooth compact connected support initial data, can be extended to [0,T*], where T* is the maximal time that the local solution u(x,t) is smooth in [0,T*), if and only if for the finite time interval [0,T*], there exist a bound M>0, so that the* ***vorticity is bounded by M, supremum norm*** *L∞ in [0,T*]:*

$$\|\omega(x,t)\|_{L_\infty} \leq M \text{ for all t in [0,T*)} \qquad \text{(eq. 18)}$$

*Then there is no maximal Blow-up time T*, and the solution exists smooth in [0,+∞)*

**Remark 4.6** Obviously if $\|\omega(x,t)\|_{L_\infty} \leq M$ , then also the integral exists and is bounded: $\int_0^{T^*} \|\omega(x,t)\|_{L_\infty} dt \leq M_1$ and the previous proposition 4.6 applies. Conversely if regularity holds, then in any interval from smoothness in a compact connected set, the vorticity is supremum bounded.

**1) PROPOSITION 4.8 (Sufficient condition for regularity)**

*The local solution u(x,t) , t in [0,T*) of the Euler or Navier-Stokes equations, with smooth compact connected support initial data, can be extended to [0,T*], where T* is the maximal time that the local solution u(x,t) is smooth in [0,T*), if for the finite time interval [0,T*], there exist a bound M>0, so that the space accelerations are bounded by M, in the supremum norm L∞ in [0,T*]:*

$$\|\nabla u(x,t)\|_{L_\infty} \leq M \text{ for all t in [0,T*)} \qquad \text{(eq. 19)}$$

*Then there is no maximal Blow-up time T*, and the solution exists smooth in [0,+∞)*

**Remark 4.7** Direct from the inequality (eq.14) and the application of the proposition 4.5. Conversely if regularity holds, then in any finite time interval from smoothness, the accelerations are supremum bounded.

**PROPOSITION 4.9 (**FEFFERMAN C. L. 2006. **Necessary and sufficient condition for regularity)**

*The local solution u(x,t) , t in [0,T*) of the Navier-Stokes equations with non-zero viscosity, and with smooth compact connected support initial data, can be extended to [0,T*], where T* is the maximal time that the local solution u(x,t) is smooth in [0,T*), if and only if*

*the velocities ||u(x,t)|| do not get unbounded as t->T*.*

*Then there is no maximal Blow-up time T*, and the solution exists smooth in [0,+∞).*

**Remark 4.8**. This is mentioned in the Official formulation of the 4th Clay Millennium problem FEFFERMAN C. L. 2006 pp.2, line 1 from below: quote "...For the Navier-Stokes equations (v>0) , if there is a solution with a finite blowup time T, then the velocities $u_i$(x,t), 1<=i<=3 become unbounded near the blowup time." The converse-negation of this is that if the velocities remain bounded near the T*, then there is no Blowup at T* and the solution is regular or global in time smooth. Conversely of course, if regularity holds, then in any finite time interval, because of the smoothness, the velocities, in a compact set are supremum bounded.

I did not find a dedicated such theorem in the books or papers that I studied, but since prof. C.L Fefferman , who wrote the official formulation of the 4th Clay Millennium problem, was careful to specify that is in the case of non-zero viscosity v>0, and not of the Euler equations as the other conditions, I assume that he is aware of a proof of it.

**PROPOSITION 4.10. (Necessary condition for regularity)**

*Let us assume that the local solution u(x,t) , t in [0,T*) of the Navier-Stokes equations with non-zero viscosity, and with smooth compact connected support initial data, can be extended to [0,T*], where T* is the maximal time that the local solution u(x,t) is smooth in [0,T*), in other words that are regular, then the trajectories-paths length l(a,t) does not get unbounded as*

*t->T*.*

**Proof:** Let us assume that the solutions is regular. Then also for all finite time intervals [0,T] , the velocities and the accelerations are bounded in the $L_\infty$, supremum norm, and this holds along all trajectory-paths too. Then also the length of the trajectories, as they are given by the formula

$$l(a_0,T) = \int_0^T \|u(x(a_0,t)\| dt \qquad \text{(eq. 20)}$$

are also bounded and finite (see e.g. APOSTOL T. 1974, theorem 6.6 p128 and theorem 6.17 p 135). Thus, if at a trajectory the lengths become unbounded as t goes to T*, then there is a blow-up. QED.

**PROPOSITION 4.11.(Physical meaningful inviscid solutions of the Euler equations are regular)**

*Let us consider the local solution u(x,t) , t in [0,T*) of the Euler equations with zero viscosity, and with smooth compact connected support initial data. If we conside,r because of zero-viscosity at every space point and at every time, as physical meaningful solutions those that also the time and space points energy dissipation density, due to viscosity, is zero or* $\|\nabla u(x,t)\|^2{}_{L_\infty}=0$ *, then , they can be extended smooth to all times [0,+∞), in other words they are regular.*

**Proof:** Direct from the PROPOSITION 4.8. QED.

**Remark 4.9.**

Similar results about the local smooth solutions, hold also for the non-homogeneous case with external forcing which is nevertheless space-time smooth of bounded accumulation in finite time intervals. Thus an alternative formulation to see that the velocities and their gradient, or in other words up to their 1st derivatives and the external forcing also up to the 1st derivatives , control the global in time existence is the next proposition. See TAO. T. 2013 Corollary 5.8

**PROPOSITION 4.12 (Maximum Cauchy development)**

*Let ($u_0$, f, T) be $H^1$ data. Then at least one of the following two statements hold:*

*1) There exists a mild $H^1$ solution (u, p, $u_0$, f, T) in [0,T] ,with the given data.*

*2)There exists a blowup time $0 < T^* < T$ and an incomplete mild $H^1$ solution*

*(u, p, $u_0$, f, $T^*$ ) up to time $T^*$ in [0, $T^*$), defined as complete on every [0,t], $t<T^*$ which blows up in the enstrophy $H^1$ norm in the sense that*

$$\lim_{t\to T^*, t<T^*} \|u(x,t)\|_{H^1_x(R^3)} = +\infty$$

**Remark 4.10** The term “almost smooth” is defined in TAO, T. 2013, before Conjecture 1.13. The only thing that almost smooth solutions lack when compared to smooth solutions is a limited amount of time differentiability at the starting time t = 0;

The term *normalized pressure*, refers to the symmetry of the Euler and Navier-Stokes equations to substitute the pressure, with another that differs at, a constant in space but variable in time measureable function. In particular normalized pressure is one that satisfies the (eq. 7) except for a measurable at a, constant in space but variable in time measureable function. It is proved in TAO, T. 2013, at Lemma 4.1, that the pressure is normalizable (exists a normalized pressure) in almost smooth finite energy solutions, for almost all times. The viscosity coefficient in these theorems of the above paper by TAO has been normalized to ν=1**.**

## 5. Conservation of the particles as a local structure of fluids in the context of continuous fluid mechanics. Proof of the regularity for fluids with conservation of particles as a local structure, and the hypotheses of the official formulation of the 4th Clay millennium problem, for the Euler and Navier-Stokes equations.

**Remark 5.1** (**Physical interpretation of the definition 5.1**) The smoothness of the particle-trajectory mapping (or displacement transformation of the points), the smoothness of the velocity field and vorticity field, is a condition that involves statements in the orders of micro scales of the fluid, larger, equal and also by far smaller than the size of material molecules, atoms and particles, from which it consists. This is something that we tend to forget in continuous mechanics, because continuous mechanics was formulated before the discovery of the existence of material atoms. On the other-hand it is traditional to involve the atoms and particles of the fluid, mainly in mathematical models of statistical mechanics. Nevertheless, we may formulate properties of material fluids in the context of continuous fluid mechanics, that reflect approximately properties and behavior in the flow of the material atoms. This is in particular the DEFINITION 5.1. **For every atom or material particle of a material fluid, we may assume around it a ball of fixed radius, called *particle range* depending on the size of the atom or particle, that covers the particle and a little bit of the electromagnetic, gravitational or quantum vacuum field around it, that their velocities and space-time accelerations are affected by the motion of the molecule or particle.** E.g. for the case water, we are speaking here for molecules of $H_2O$, that are estimated to have a diameter of 2.75 angstroms or 2r= 2.75*10^(-10) meters, we may define as water molecule **particle range** the balls $B(r_0)$ of radius $r_0$=3*10^(-10) meters around the water molecule. As the fluid flows, especially in our case here of incompressible fluids, the shape and size of the molecules do not change much, neither there are significant differences of the velocities and space-time accelerations of parts of the molecule. Bounds $\delta_u$ $\delta_\omega$ of such differences remain constant as the fluid flows. **We may call this effect as the *principle of conservation of particles* as a local structure.** This principle must be posed in equal setting as the energy conservation and incompressibility together with the Navier-Stokes or Euler equations. Of course, if the fluid is say of solar plasma matter, such a description would not apply. Nevertheless, then incompressibility is hardly a property of it. But if we are talking about incompressible fluids that the molecule is conserved as well as the atoms and do not change atomic number (as e.g. in fusion or fission) then this principle is physically valid. The principle of conservation of particles as a local structure, blocks the self-similarity effects of concentrating the energy and turbulence in very small areas and creating thus a Blow-up. It is the missing invariant in the discussion of many researchers about superctitical, critical and subcritical invariants in scale transformations of the solutions.

The next DEFINITION 5.1 formulates precisely mathematically this principle for the case of incompressible fluids.

**DEFINITION 5.1. (Conservation of particles as local structure in a fluid)**

*Let a smooth solution of the Euler or Navier-Stokes equations for incompressible fluids, that exists in the time interval [0,T). We may assume initial data on all of $R^3$ or only on a connected compact support $V_0$. For simplicity let us concentrate only on the latter simpler case. Let us denote by F the displacement transformation of the flow Let us also denote by g the partial derivatives of 1st order in space and time , that is $\left|\partial_x^a \partial_t^b u\ (x)\right|$, $|\alpha|=1$, $|b|<=1$,and call then space-time accelerations . We say that there is **conservation of the particles** **in the interval [0,T) in** a derivatives homogenous setting, as a local structure of the solution if and only if:*

*There is a small radius r, and small constants $\delta_x$ , $\delta_u$ , $\delta_\omega$ , >0 so that for all t in [0,T) there is a finite cover $C_t$ (in the case of initial data on $R^3$ · it is infinite cover, but finite on any compact subset) of $V_t$ , from balls B(r) of radius r, called* ***ranges of the particles*** *, such that:*

*4) For an $x_1$ and $x_2$ in a ball B(r) of $V_s$ , s in [0,T), $||F(x_1)-F(x_2)||<=r+ \delta_x$ for all t>=s in [0,T).*

*5) For an $x_1$ and $x_2$ in a ball B(r) of $V_s$ ,s in [0,T), $||u(F(x_1))-u(F(x_2))||<= \delta_u$ for all t >=s in [0,T).*

*6) For an $x_1$ and $x_2$ in a ball B(r) of $V_s$ , s in [0,T), $||g(F(x_1))-g(F(x_2))||<= \delta_\omega$ for all t >=s in [0,T).*

*If we state the same conditions 1) 2) 3) for all times t in [0,+∞) , then we say that we have the* ***strong version*** *of the conservation of particles as local structure.*

**PROPOSITION 5.1 (Velocities on trajectories in finite time intervals with finite total variation, and bounded in the supremun norm uniformly in time.)**

*Let $u_t$ : V(t) -> $R^3$ be smooth local in time in [0,T*) ,velocity fields solutions of the Navier-Stokes or Euler equations, with compact connected support V(0) initial data, finite initial energy E(0) and conservation of particles in [0,T*) as a local structure . The [0,T*) is the maximal interval that the solutions are smooth. Then for t in [0,T*) and x in V(t), the velocities are uniformly in time bounded in the supremum norm by a bound M independent of time t.*

$$\|u(x,t)\|_{L_\infty} = \sup_{x\in V(t)} \|u(x,t)\| \leq M \text{ for all } t \text{ in } [0,T^*).$$

*Therefore, the velocities on the trajectory paths, in finite time intervals are of bounded variation and the trajectories in finite time interval, have finite length.*

**1st Proof (Only for the Navier-Stokes Equations):** Let us assume, that the velocities are unbounded in the supremum norm, as t converges to T*. Then there is a sequence of times $t_n$ with $t_n$ converging to time T* , and sequence of corresponding points $x_n$ ($t_n$ ), for which the norms of the velocities $||u( x_n (t_n ), t_n)||$ converge to infinite.

$$\lim_{n\to+\infty} || u(x(x_n,t_n),t_n) ||= +\infty . \quad \text{(eq.21)}$$

From the hypothesis of the conservation of particles as a local structure of the smooth solution in [0,T*), for every $t_n$ There is a finite cover $C_{tn}$ of particle ranges, of $V_{tn}$ so that $x_n$ ($t_n$ ) belongs to one such ball or particle-range $B_n(r)$ and for any other point $y(t_n)$ of $B_n(r)$, it holds that $||u(x_n (t_n ), t_n)-u(y(t_n),t_n)||<= \delta_u$ . Therefore

$$||u(x_n (t_n ), t_n)||- \delta_u <= ||u(y(t_n),t_n)||<= ||u(x_n (t_n ), t_n)||+ \delta_u \quad \text{(eq.22)}$$

for all times $t_n$ in [0,T*) .

By integrating spatially on the ball $B_n(r)$, and taking the limit as n->+∞ we deduce that

$$\underset{n\to+\infty}{Lim}\int_{B_n}\|u\|dx=+\infty$$

But this also means as we realize easily, that also

$$\underset{n\to+\infty}{Lim}\int_{B_n}\|u\|^2 dx=+\infty \qquad \text{(eq. 23)}$$

Which nevertheless means that the total kinetic energy of this small, but finite and of constant radius, ball, converges to infinite, as $t_n$ converges to T*. This is impossible by the finiteness of the initial energy, and the conservation of energy. Therefore the velocities are bounded uniformly ,in the supremum norm, in the time interval [0,T*).

Therefore the velocities on the trajectory paths, are also bounded in the supremum norm , uniformly in the time interval [0,T*). But this means by PROPOSITION 4.9 that the local smooth solution is regular , and globally in time smooth, which from PROPOSITION 4.8 means that the Jacobian of the 1st order derivatives of the velocities are also bounded in the supremum norm uniformly in time bounded in [0,T*). Which in its turn gives that the velocities are of bounded variation on the trajectory paths (see e.g. APOSTOL T. 1974 , theorem 6.6 p128 and theorem 6.17 p 135) and that the trajectories in have also finite length in [0,T*), because the trajectory length is given by the formula $l(a_0,T)=\int_0^T\|u(x(a_0,t)\|dt$ .

QED.

**2nd Proof (Both for the Euler and Navier-Stokes equations):** Instead of utilizing the condition 2) of the definition 5.1, we may utilize the condition 3). And we start assuming that the Jacobian of the velocities is unbounded in the supremum norm (instead of the velocities), as time goes to the Blow-up time T*. Similarly we conclude that the energy dissipation density at a time on balls that are particle-ranges goes to infinite, giving the same for the total accumulative in time energy dissipation (see (eq. 11), which again is impossible from the finiteness of the initial energy and energy conservation. Then by PROPOSITION 4.8 we conclude that the solution is regular, and thus also that the velocities are bounded in the supremum norm, in all finite time intervals. Again, we deduce in the same way, that the total variation of the velocities is finite in finite time intervals and so are the lengths of the trajectories too.

QED.

**PROPOSITION 5.2 (Global regularity as in the 4th Clay Millennium problem).**

*Let the Navier-Stokes or Euler equations with smooth compact connected initial data, finite initial energy and conservation of particles as local structure. Then the unique local in time solutions are also regular (are smooth globally in time).*

**Proof:** We apply the PROPSOITION 5.1 above and the necessary and sufficient condition for regularity in PROPSOITION 4.9 (which is only for the Navier-Stokes equations). Furthermore, we apply the part of the 2d proof of the PROPOSITION 5.1, that concludes regularity from PROPSOITION 4.8 which holds for both the Euler and Navier-Stokes equations.

QED.

## 6. Bounds of measures of the turbulence from length of the trajectory paths, and the total variation of the velocities, space acceleration and vorticity. The concept of homogeneous smoothness.

**Remark 6.1** In the next we define **a measure of the turbulence** of the trajectories, of the velocities, of space-time accelerations and of the vorticity, through the **total variation** of the component functions in finite time intervals. This is in the context of deterministic fluid dynamics and not stochastic fluid dynamics. We remark that in the case of a blowup the measures of turbulence below will become infinite.

**DEFINITION 6.1 (The variation measure of turbulence)**

Let smooth local in time in [0,T] solutions of the Euler or Navier-Stokes equations. The total length L(P) of a trajectory path P, in the time interval [0,T] is defined as **the variation measure of turbulence of the displacements** on the trajectory P, in [0,T].The total variation TV(||u||) of the norm of the velocity ||u|| on the trajectory P in [0,T] is defined as **the variation measure of turbulence of the velocity** on the trajectory P in [0,T]. The total variation TV(g) of the space-accelerations g (as in DEFINITION 5.1) on the trajectory P in [0,T] is defined as **the variation measure of turbulence of the space-time accelerations** on the trajectory P in [0,T]. The total variation TV(||ω||) of the norm of the vorticity ||ω|| on the

trajectory P in [0,T] is defined as **the variation measure of turbulence of the vorticity** on the trajectory P in [0,T].

**PROPOSITION 6.1 Conservation in time of the boundedness of the maximum turbulence, that depend only on the initial data and time lapsed.**

*Let the Euler or Navier-Stokes equations with smooth compact connected initial data finite initial energy and conservation of the particles as a local structure. Then for all times t, there are bounds $M_1(t)$, $M_2(t)$, $M_3(t)$, so that the maximum turbulence of the trajectory paths, of the velocities and of the space accelerations are bounded respectively by the above universal bounds, that depend only on the initial data and the time lapsed.*

**Proof:** From the PROPOSITIONS 5.1, 5.2 we deduce that the local in time smooth solutions are smooth for all times as they are regular. Then in any time interval [0,T], the solutions are smooth, and thus from the PROPOSITION 4.8, the space acceleration g, are bounded in [0,T], thus also as smooth functions their total variation TV(g) is finite, and bounded. (see e.g. APOSTOL T. 1974, theorem 6.6 p128 and theorem 6.17 p 135).From the PROPOSITION 4.7, the vorticity is smooth and bounded in [0,T], thus also as smooth bounded functions its total variation TV(||ω||) is finite, and bounded on the trajectories. From the PROPOSITION 4.9, the velocity is smooth and bounded in [0,T], thus also as smooth bounded functions its total variation TV(||u||) is finite, and bounded on the trajectories. From the PROPOSITION 4.10, the motion on trajectories is smooth and bounded in [0,T], thus also as smooth bounded functions its total variation which is the length of the trajectory path L(P) is finite, and bounded in [0,T].In the previous theorems the bounds that we may denote them here by $M_1(t)$, $M_2(t)$, $M_3(t)$, respectively as in the statement of the current theorem, depend on the initial data, and the time interval [0,T]. QED.

**Remark 6.2. (Homogeneity of smoothness relative to a property P.)** There are many researchers that they consider that the local smooth solutions of the Euler or Navier-Stokes equations with smooth Schwartz initial data and finite initial energy, (even without the hypothesis of conservation of particles as a local structure) are general smooth functions. But it is not so! They are special smooth functions with the remarkable property that there are some critical properties $P_i$ that if such a property holds in the time interval [0,T) for the coordinate partial space-derivatives of 0, 1, or 2 order , then this property holds also for the other two orders of derivatives. In other words if it holds for the 2 order then it holds for the orders 0, 1 in [0,T) . If it holds for the order 1, then it holds for the orders 0, 2 in [0,T]. If it holds for the order 0, them it holds also for the orders 1,2 in [0,T]. This pattern e.g. can be observed for the property $P_1$ of uniform boundedness in the supremum norm, in the interval [0,T*) in the PROPOSITIONS 4.5-4.10 . But one might to try to prove it also for a second property $P_2$ which is the **finitness of the total variation** of the coordinates of the partial derivatives, or even other properties P3 like **local in time Lipchitz conditions**. This creates a strong bond or coherence among the derivatives and might be called ***homogeneous smoothness relative to a property P.*** We may also notice that the formulation of the conservation of particles as local structure is in such a way, that as a property, it shows the same pattern of homogeneity of smoothness relative to the property of uniform in time bounds $P_4$ ,1), 2), 3) in the DEFINITION 5.1. It seem to me though that even this strong type of smoothness is not enough to derive the regularity, unless the homogeneity of smoothness is relative to the property $P_4$ , in other words the conservation of particles as a local structure.

## 7. Epilogue

I believe that the main reasons of the failure so far in proving of the 3D global regularity of incompressible flows, with reasonably smooth initial conditions like smooth Schwartz initial data, and finite initial energy, is hidden in the difference of the physical reality ontology that is closer to the ontology of statistical mechanics models and the ontology of the mathematical models of continuous fluid dynamics.

Although energy and momentum conservation and finiteness of the initial energy are easy to formulate in both types of models, the conservation of particles as type and size is traditionally formulated only in the context of statistical mechanics. By succeeding in formulating approximately in the context of the ontology of continuous fluid mechanics the conservation of particles during the flow, as local structure, we result in being able to prove the regularity in the case of 3 dimensions which is what most mathematicians were hoping that it holds.

## CHAPTER 6

# PROPOSED SOLUTION OF THE MILLENNIUM PROBLEM ON THE NAVIER-STOKES EQUATIONS BASED ON THE PHYSICAL PRINCIPLE OF PARTICLES STRUCTURE CONSERVATION OF THE FLUIDS.

**Konstantinos E. Kyritsis***

* Associate Prof. of University of Ioannina, Greece. Dept. Accounting-Finance Psathaki Preveza 48100

Corresponding Author: ckiritsi@uoi.gr, C_kyritsis@yahoo.com

**Abstract** In this article we propose a solution of the Millennium problem about the Navier-Stokes equations. The solution is in favor of the regularity (no finite time blowup) for the homogeneous case (no external forces) under the additional physical hypothesis of volume preserving particle structure conservation, which is formulated in the context of fluid dynamics. Strictly speaking it is not the solution of the original mathematical problem as it is added the particle structure conservation. Still it is the closest possible solution from the point of view of physics and expectations of the physicists based on our current understanding of matter, not as infinite divisible but as particle based even for fluid dynamics.

---

## 1. Introduction

The Clay millennium problem about the Navier-Stokes equations is one of the 7 famous problem of mathematics that the Clay Mathematical Institute has set a high monetary award for its solution. It is considered a difficult problem as it has resisted solving it for almost a whole century. The Navier-Stokes equations are the equations that are considered to govern the flow of fluids, and had been formulated long ago in mathematical physics before it was known that matter consists from atoms. So actually they formulate the old infinite divisible material fluids. Although it is known that under its assumptions of the millennium problem the Navier-Stokes equations have a unique smooth and local in time solution, it was not known if this solution can be extended smoothly and globally for all times, which would be called the **regularity of the Navier-Stokes equations in 3 dimensions.** The corresponding case of regularity in 2 dimensions has long ago been proved to hold but the 3-dimensions had resisted proving it. There mainly two groups of people with different expectations about this millennium problem. At first it is the group of mathematical physicists, who consider it natural that there should not exist a blow-up in finite time in 3 dimesions either. Then there is the group of mathematicians of partial differential equations, who think that it is very probable that it might exist a blow up in finite time. Both groups seem to be right! The reason is the next. There are physical aspects of the Millennium Problem, that are not formulated mathematically, in the setting of the problem. For example, the energy conservation of the heat of the fluid, is not included in the formulation of the problem at all. There is no symbol or magnitude for the temperature in the formulation of the problem. Still a physicist thinks, that the conservation of the energy of the heat is something that does hold and can use it in the arguments. Nevertheless, a mathematician, considers the millennium problem, as a mathematical problem of partial differential equations, and considers as holding, only what is mathematically stated in the formulation. I belong to the group of mathematical physicists, and that is how I propsoe a solution of the millennium problem here. In particular, I use the **physical principle of particles structure, conservation**, in the context of fluid mechanics.

This principle must be posed in equal setting as the energy conservation and incompressibility together with the Navier-Stokes or Euler equations. Of course if the fluid is say of solar plasma matter, such a description would not apply. Nevertheless then incompressibility is hardly a property of it. But if we are talking about incompressible fluids that the molecule is conserved as well as the atoms and do not change atomic number (as e.g. in fusion or fission) then this principle is physically valid. The principle of conservation of particles as a local structure, blocks the self-similarity effects of concentrating the energy and turbulence in very small areas and creating thus a Blow-up. It is the missing invariant in the discussion of many researchers about superctitical, critical and subcritical invariants in scale transformations of the solutions.

**2. The formulation of the millennium problem and the 4 sub-problems (A) , (B), (C), (D)**

In this paragraph we highlight the basic parts of the standard formulation of the 4th Clay millennium problem as in [1] Fefferman C.L. 2006.

The **Navier-Stokes** equations are given by (by R we denote the field of the real numbers, ν>0 is the density normalized viscosity coefficient )

$$\frac{\partial}{\partial t}u_i + \sum_{j=1}^{n} u_j \frac{\partial u_i}{\partial x_j} = -\frac{\partial p}{\partial x_i} + \nu \Delta u_i \qquad (x \in R^3 , t>=0 , n=3) \qquad \text{(eq.2.1 )}$$

$$div u = \sum_{i=1}^{n} \frac{\partial u_i}{\partial x_i} = 0 \qquad (x \in R^3 , t>=0 , n=3) \qquad \text{(eq. 2.2 )}$$

with initial conditions $u(x,0)=u^0(x)$ $\quad$ **x∈R³**

and $u^0$ (x) ∈ $C^\infty$ divergence-free vector field on **R³** (eq. 2.3)

If ν=0 then we are taking about the **Euler** equations and **inviscid** case.

$\Delta = \sum_{i=1}^{n} \frac{\partial^2}{\partial x_i^2}$ is the Laplacian operator .The Euler equations are (eq2.1), (eq 2.2), (eq 2.3) when ν=0.

It is reminded to the reader, that in the equations of Navier-Stokes, as in (eq. 2.1) the density ρ, is constant, it is custom to normalized to 1 and omit it.

For physically meaningful solutions we want to make sure that $u^0(x)$ does not grow large as $|x| \to \infty$. This is set by defining $u^0(x)$ , and f(x,t) and called in this paper **(smooth) Schwartz initial conditions** , to satisfy

$\left|\partial_x^a u^0(x)\right| \le C_{a,K}(1+|x|)^{-K}$ on **R³** for any *α* and K (eq.2.4 )

(Schwartz used such functions to define the space of Schwartz distributions)

**Remark 2.1.** It is important to realize that smooth Schwartz initial velocities after

(eq 2.4) will give that the initial vorticity $\omega_0$ =curl($u^0$) , in its supremum norm, is bounded over all 3-space.

We also assume **external forces** that are Schwartz in space and time:

$\left|\partial_x^a \partial_t^m f \quad (x,t)\right| \le C_{a,m,K}(1+|x|+t)^{-K}$on **R³**× $[\mathbf{0},+\infty)$for any *α,m,*K (eq.2.5 )

We accept as physical meaningful solutions only if it satisfies

p, u $\in C^\infty(R^3 \times [0,\infty))$ (eq.2.6 )

and

$\int_{\Re^3} |u(x,t)|^2 dx < C$ for all t>=0 (**Bounded or finite energy**) (eq.2.7 )

**Lemma 2.0.** *from the smooth Schwartz initial data of the formulation of the 4th Clay Millennium problem we derive that the supremum of the vorticity at t=0, symbolized by $F_\omega$ exist and it is finite.*

**Proof.** We use the remark 2.1 and of course the inequality |a-b|<=|a|+|b|, at the definition of the vorticity by the partial derivatives. QED.

**Remark 2.2** It is important to realize that smooth external force (densities) with the Schwartz property as in (eq.2.5) , have not only a rule for upper bounded spatial partial derivatives but also the same rule for time upper bounded partial derivatives.

**Remark 2.3** We must stress here that imposing smoothness of the coordinate functions of velocities and external forces of the initial t=0 data and later time t data in the Cartesian coordinates plus and Schwartz condition as in (eq 2.5) is not equivalent with imposing similar such smoothness of the coordinate functions and conditions in the cylindrical or spherical coordinates. We will give in the paragraph 4, remark 4.5 an example of a strange blowup, where at any time t>0 the coordinates of the velocities are smooth and bounded in all space as functions in the polar coordinates and still the vorticity has infinite singularity at zero.

Alternatively, to rule out problems at infinity, we may look for spatially periodic solutions of (2.1), (2.2), (2.3). Thus we assume that $u^0(x)$ , and f(x,t) satisfy

$u^0(x+e_j)= u^0(x)$, $f(x+e_j,t)= f(x,t)$, $p(x+e_j,0)=p(x,0)$, for 1<=j<=3 (eq.2.8 )

($e_j$ is the jth unit vector in $R^3$)

In place of (2.4) and (2.5), we assume that $u^0(x)$, is smooth and that

$\left|\partial_x^a \partial_t^m f \quad (x,t)\right| \leq C_{a,m,K}(1+t)^{-K}$ on $\mathbf{R^3} \times [\mathbf{0}, +\infty)$ for any $\alpha, m,$ K (eq.2.9 )

We then accept a solution of (2.1), (2.2) , (2.3) as physically reasonable if it satisfies

$u(x+e_j,t)= u(x, t)$, $p(x+e_j, t)=p(x,t)$, on $\mathbf{R^3} \times [\mathbf{0}, +\infty)$ for 1<=j<=3 (eq.2.10 )

and p, u $\in C^\infty(R^3 \times [0,\infty))$ (eq.2.11 )

In the next paragraphs we may also write $u_0$ instead of $u^0$ for the initial data velocity.

We denote Euclidean balls by $B(a,r) := \{x \in R^3 : ||x-a|| \leq r\}$, where ||x|| is the Euclidean norm.

**The 4 sub-problems or conjectures of the millennium problem are the next:**

**(Conjecture A) Existence and smoothness of Navier-Stokes solution on $R^3$.**

*Take v>0 and n=3. Let $u_0(x)$ be any smooth, divergent-free vector field satisfying (2.4). Take f(x,t) to be identically zero. Then there exist smooth functions p(x,t) , u(x,t) on $R^3x[0,+\infty)$ that satisfy (2.1), (2.2), (2.3) , (2.6) , (2.7).*

**(Conjecture B) Existence and smoothness of Navier-Stokes solution on $R^3/Z^3$.**

*Take v>0 and n=3. Let $u_0(x)$ be any smooth, divergent-free vector field satisfying (2.8); we take f(x,t) to be identically zero. Then there exist smooth functions p(x,t) , u(x,t) on $R^3x[0,+\infty)$ that satisfy (2.1), (2.2), (2.3) , (2.10) , (2.11).*

**(Conjecture C) Breakdown of Navier-Stokes solution on $R^3$**

*Take v>0 and n=3. Then there exist a smooth, divergent-free vector field $u_0(x)$ on $R^3$ and a smooth f(x,t) on $R^3x[0,+\infty)$ satisfying (2.4), (2.5) for which there exist no smooth solution (p(x,t) ,u(x,t)) of (2.1), (2.2), (2.3) , (2.6) , (2.7) on $R^3x[0,+\infty)$.*

**(Conjecture D) Breakdown of Navier-Stokes solution on $R^3/Z^3$**

*Take v>0 and n=3. Then there exist a smooth, divergent-free vector field $u_0(x)$ on $R^3$ and a smooth f(x,t) on $R^3x[0,+\infty)$ satisfying (2.8), (2.9) for which there exist no smooth solution (p(x,t) ,u(x,t)) of (2.1), (2.2), (2.3) , (2.10) , (2.11) on $R^3x[0,+\infty)$.*

**Remark 2.4.** It is stated in the same formal formulation of the Clay millennium problem by C. L. Fefferman see [1] Fefferman C.L. 2006 (see page 2nd line 5 from below) that the conjecture (A) has been proved to holds locally. "..if the time internal $[0,\infty)$, is replaced by a small time interval [0,T), with T depending on the initial data....". In other words there is $\infty>T>0$, such that there exists a unique and smooth solution $u(x,t) \in C^{\infty}(R^3 \times [0,T))$. See also [3] A.J. Majda-A.L. Bertozzi ,Theorem 3.4 pp 104. In this paper, as it is standard almost everywhere, the term smooth refers to the space $C^{\infty}$

In the next the $|| \, ||_m$ is the corresponding Sobolev spaces norm and . We denote by $V^m$ ={u in $H^m(R^n)$ and divu=0} where $H^m(R^n)$ are the Sobolev spaces with the $L^2$ norm.

We must mention that in A.J. Majda-A.L. Bertozzi [3] ,Theorem 3.4 pp 104, Local in Time existence of Solutions to the Euler and Navier-Stokes equations it is proved that indeed if the initial velocities belong to $V^m$ m>=[3/2]+2 there exist unique smooth solutions locally in time [0,t]. Here, in the formulation of the millennium problem the hypotheses of smooth with Schwartz condition initial velocities satisfies this condition therefore we have the existence and uniqueness of smooth solution locally in time, both in the non-periodic and the periodic setting without external forcing (homogeneous case).

The existence and uniqueness of a smooth solutions locally in time is stated in the formulation by C.L. Fefferman [1] for the homogeneous cases and conjectures (A), (B). When a smooth Schwartz condition external force is added (inhomogeneous case) , it is natural to expect that also there should exist a local in time unique sooth solution. **And this is indeed provable with the existing bibliography** (The standard theorems are due to **Fujita–Kato**, **Kato**, **Leray**, etc and presented in many textbooks such as Temam [6], Ladyzhenskaya [7] , Majda [3] etc ).

**Lemma 2.1** *The smooth Schwartz initial data and the rest of the hypotheses in the formulation of the $4^{th}$ Cay Millennium problem for the homogeneous case and conjectures (A), (B), and non-homogeneous case and conjectures (C), (D), guarantee, that it exists a unique smooth solution locally in time [0,t) , for the Navier-Stokes equations (without or with external forces).*

We state here also two, very well-known criteria of no blow-up and regularity.

In the next theorem the $|| \ ||_m$ is the corresponding Sobolev spaces norm and . We denote by $V^m$ ={u in $H^m(R^n)$ and divu=0} where $H^m(R^n)$ are the Sobolev spaces with the $L^2$ norm.

**Theorem 2.1 Velocities Sobolev norm sufficient condition of regularity.** *Given an initial condition $u_0 \in V^m$ m>=[3/2]+2=3.5 e.g. m=4 , then for any viscosity v>=0 . there exists a maximal time T* (possibly infinite) of existence of a unique smooth solution $u \in C([0,T^*];V^m) \cap C^1([0,T^*];V^{m-2})$ to the Euler or Navier-Stokes equation. Moreover, if T*<+∞ then necessarily $\lim_{t\to T^*} ||u(.,t)||_m = +\infty$.*

**Proof:** See A.J. Majda-A.L. Bertozzi [3], Corollary 3.2 pp 112). QED

**Remark 2.5** Obviously this proposition covers the periodic case too.

**Theorem 2.2 (Bealo-Kato-Majda) Supremum of vorticity sufficient condition of regularity**

*Let the initial velocity $u_0 \in V^m$ m>=[3/2]+2 , e.g. m=4, so that there exists a classical solution $u \in C^1([0,T]; C^2 \cap V^m)$ to the 3D Euler or Navier-Stokes equations. Then :*

*(i) If for any T>0 there is $M_1$ >0 such that the vorticity ω=curl(u) satisfies*

$$\int_0^T |\omega(.,\tau)|_{L^\infty}\, d\tau \leq M_1$$

*Then the solution u exists globally in time, $u \in C^1([0,+\infty]; C^2 \cap V^m)$*

*(ii) If the maximal time $T^*$ of the existence of the solution $u \in C^1([0,T]\,;\, C^2 \cap V^m)$ is finite, then necessarily the vorticity accumulates so rapidly that*

$$lim_{t \to T*} \int_0^T |\omega(.,\tau)|_{L^\infty}\, d\tau = +\infty \qquad \textit{(eq. 2.12)}$$

**Proof:** See A.J. Majda-A.L. Bertozzi [3], Theorem 3.6 pp 115, $L^\infty$ vorticity control of regularity. QED.

**Remark 2.6** Obviously this proposition covers the periodic case too. But furthermore it is important to know, that the above Bealo-Kato-Majda, holds not only for the **Navier-Stokes** but the for the **Euler** equations and furthermore **without or with external forcing under the Schwartz boundedness conditions above.**

## 3. What is that maybe we do not understand very well with the Navier-Stokes equations and their physics?

In statistical mechanical models of incompressible flow, we have the realistic advantage of finite many particles, e.g. like balls B(x,r) with finite diameter r. These particles as they flow in time, remain particles of the same nature and size and the velocities and inside them remain approximately constant.

F. g. for the case water, we are speaking here for molecules of $H_2O$, that are estimated to have a diameter of 2.75 angstroms or 2r= 2.75*10^(-10) meters. The Representative Elementary Volume or REV is usually assumed about **$10^{-15}$ meters.** But if we go deaper e.g. for electrons, then the size grows smaller between **$10^{-23}$ meters** and $10^{-15}$ meters. Furthermore all the electrons , protons and neutrons are assumed to spin with about the same frequency which is about 1 terahertz. Or as **vorticity** ω it would be 2π times that or **ω=6.26 Terraherz.** This frequency is what makes our molecular matter physical reality coherent, and it make increase or decrease simultaneously in large areas from deeper factors, that might as well be related to human consciousness too, as some experiments prove.

Fluid dynamics from the strictly mathematical point of view goes smaller than any scale, as it was conceived at the centuries that matter was believed to be infinite divisible. It is somehow a shame that we still use such type of models! In practice when we apply e.g. Navier-Stokes fluid dynamics in water fluids we assume that we may go as small as the Representative Elementary Volume or REV is usually assumed about **$10^{-15}$ meters.** In the principle of the particle structure conservation that we define and use in this article, we assume that the vorticity ω of the **Representative Elementary Volumes but the point vorticities too are uniformly modulated and influences by the constant spin of ω=6.26 Terraherz of the electrons, neutrons and protons.**

Because space and time dimensions in classical fluid dynamics goes in orders of smallness , smaller and at least as small as the real physical molecules , atoms and particles of the fluids, this might suggest imposing too, such conditions resembling uniform continuity conditions. In the case of continuous fluid dynamics models such natural conditions, emerging from the particle nature of material fluids, together with the energy conservation, the incompressibility and the momentum conservation, as laws conserved in time, may derive the regularity of the local smooth solutions of the Navier-Stokes equations.

In [8] it is defined **a concept of particle structure and conservation of particle structure** in incomprehensible fluids, which is used to propose a solution of the Millennium problem about the Navier-Stokes eqautions in favor of the regularity (no blow-up). But here we define a simple concept of particle structure and conservation of particle structure, based only on the **vorticity ,** which gives as simpler solution of the millennium problem about the Navier-Stokes equations in favor of the regularity.

**Definition 3.0 The physical principle of particle structure conservation**

*Let a non-compressble fluid, which satisfies the Euler or Navier-Stokes equations, and flows smoothly on all of 3D space and in the time interval [0, T). it may be defined as non-periodic or periodic flow, and without forces. We say that it has at time t=0* ***particle structure*** *relative to vorticity* ***at smallness scale ε*** *iff for any point x, where it has vorticity ω(x,0), there is a diffeomorphic image B(0) of a spherical ball, with volume equal to |B(0)|=ε , which contains x , and the 3D volume average of the vorticity on B(0) is equal to ω(x,0). In symbols*

$$\int_{B(0)} \omega dv = \omega(x,0)$$

*We say that we have* ***volume preserving particle structure conservation*** *during [0,T) , if the same happens with the* ***same ε,*** *for all later times t in [0,T).*

We may define an alternative and more weak condition to reflect the particle structure conservation of a classical Fluid, which we call Convolution Weak particle structure. We state it here we can derive regularity with this too, but we will not as we will use in the next only the simple geometric version **particle structure** relative to vorticity **at smallness scale ε.**

**Definition 3.0.1 Convolution Weak particle structure.**

We say that a fluid has **weak convolution particle structure** if and only if there is a fixed length ε >0, which represents the scale where the particles influence the fluid, such that

$$\sup_{t\in[0,T)} (\sup_{x\in R^3} | \int_{R^3} \omega(y,t) G_\varepsilon(x-y)dy |) < +\infty$$

Where $G_\varepsilon = \varepsilon^{-3} G(\frac{x}{\varepsilon}), \int_{R^3} G(x)dx = 1$ ,G>=0 smooth rapidly decaying at infinite,

So that $G_\varepsilon$ is a **smooth averaging kernel for convolution at length scale ε,** giving a **local spatial average of vorticity over a ball of radius about ε**.

Obviously here we may take as ε the volume of a Representative Elementary Volume or REV is usually assumed about **$10^{-15}$ meters,** thus ($10^{-15}$ meters)$^3$=ε or

**ε=10^(-45) $m^3$.**

Such a physical principle may seem very strong and unusual, to the standards of fluid dynamics, that pretends that the matter is infiniote divisible, but it has physical anchors. To see that **it is compatible with the Navier-Stokes equations** , we only need to notice that the irrotational flows (or potential flows) do satisfy it, and obviously from the physical point of view, it must be applied for a smallness scale ε, much larger than the ε=10^(-45) $m^3$.

It is also important to realize that , we define volume preserving particle structure conservation, thus **we must have non-compressible volume preserving flows. For compressible fluids this definition does not apply.**

Besides the forgotten conservation law of finite particles, there are **two more forgotten laws of conservation or invariants**. The first of them is the obvious that during the flow, the physical measuring units dimensions (dimensional analysis) of the involved physical quantities (mass density, velocity, vorticity, momentum, energy, force point density, pressure, etc.) are conserved. It is not very wise to eliminate the physical magnitudes interpretation and their dimensional analysis when trying to solve the millennium problem, because the dimensional analysis is a very simple and powerful interlink of the involved quantities and leads with the physical interpretation, to a transcendental shortcut to symbolic calculations. **By eliminating the dimensional analysis, we lose part of the map to reach our goal.**

The 3rd forgotten conservation law or invariant, is related to the viscosity (friction). Because we do know that at each point (pointwise), the viscosity is only subtracting kinetic energy, with an irreversible way, and converting it to thermal energy, (negative energy point density), and this is preserved in the flow, (it can never convert thermal energy to macroscopic kinetic energy), we know that its sign does not change too, it is a flow invariant, so the integrated 1D or 2D work density is always of the same sign (negative) and as sign, an invariant of the flow. The **conservation or invariance of the sign of work density by the viscosity (friction)** is summarized in the lemma 3.1 below.

**Viscosity irreversibility Principle 3.2 The absolute physical principle of irreversibility of the dissipated kinetic energy by viscosity, for non-compressible fluids as forgotten invariant.**

*For incompressible viscous fluids, the dissipated energy by viscosity, only subtracts from the kinetic energy and convert it to thermal. It is not reversible, it cannot add to the kinetic energy. And this applies not only to the total 3D space and time integral, but also to the 2D and 1D densities of rate of energy dissipation (absolute irreversibility=both global as 3D integral but also local as 2D, 1D or point values densities of rates) . As a consequence, on a trajectory path, although the viscosity force density is not aligned with the velocity, its projection on the trajectory path, always opposes the pressure force densities projection on the trajectory path. Similarly, at the circulation of velocities multiplied by the constant density (as a 1D circular momentum density), the work density of the viscosity force densities on the same circle, will always subtract kinetic energy and reduce the momentum density. If the mathematics of solutions of the NSE equations, allow violations of the above irreversibility of the dissipation energy, then they are not acceptable as physical meaningful solutions. In the same way that solutions where the kinetic energy is increasing are not acceptable as physically meaningful solutions (as L. Fefferman remarks in [3] in the formulation of the millennium problem).*

Although we shall use the physical principle of volume preserving particle structure conservation, we shall not use in the next in this article the physical principle of irreversibility of the dissipated kinetic energy.

**4. The particle structure conservation solves that Millennium problem about the Navier-Stokes equations in favor of the regularity (no-blowup).**

We prove a basic theorem 4.1 in this paragraph.

We remind the reader, that the hypotheses of the millennium problem, involve only **finite kinetic energy**. Thus for physical reasons **the energy dissipation is also finite** , and thus the integral which calculates it through the vorticity at time t>0 (enstrophy) and as far as the local in time unique smooth solution exists, It exists also and it is finite

$$(\int_0^t (\int_{R^3} ||\omega(x,s)||^2 \, dx)ds) < +\infty \qquad \textit{(eq 4.0)}$$

**Theorem 4.1** *Under the hypotheses as in the conjectures A and B of the formulation of the millennium problem (homogeneous case without external forces) and the existence of a unique smooth solution in the time interval [0,T*) and the additional hypothesis that the fluid has volume preserving particle structure conservation as in the definition 3.1, at smallness scale ε, then the finite accumulation of the supremum norm of the vorticity in the time interval [0,T) is bounded by a constant multiple of the finite accumulation of the enstrophy through the L2 norm , which it is finally finite.*

*If by the lower right indices of the symbol ||,|| we denote if it is a L2 or a L-infinite (supremum) norm then in symbols*

$$||\omega(t)||_{L_\infty(R^3)} \le \varepsilon^{(-1/2)} ||\omega(t)||_{L^2(R^3)} \qquad \textit{(eq 4.1)}$$

*Thus*

$$\int_0^T ||\omega(.,\tau)||_{L^\infty} \, d\tau \le \varepsilon^{(-1/2)} \int_0^t ((\int_{R^3} ||\omega(s,\tau)||^2 \, ds)^{(1/2)}) d\tau = \varepsilon^{(-1/2)} t^{(1/2)} (\int_0^t (\int_{R^3} ||\omega(x,s)||^2 \, dx)ds) < +\infty$$

*(eq. 4.2)*

**Proof:**

We fix an s <T* and a x in $R^{3,}$ in 3D space. By the particle structure conservation

we have that ω(x,s)=( 1/ε) $\int\limits_{B(x,s)} \omega(y,s)dy$ After taking the Euclidean norms

$\|\omega(x,s)\| \le \frac{1}{\varepsilon} \int\limits_{B(x,s)} \|\omega(y,s)\| dy$ Then by the Cauchy-Schwartz inequality we get

$$\int\limits_{B(x,s)} \|\omega\| dy \le | B(x,s) |^{(1/2)} ( \int\limits_{B(x,s)} \|\omega\|^2 dy)^{(1/2)} = \varepsilon^{(1/2)} ( \int\limits_{B(x,s)} \|\omega\|^2 dy)^{(1/2)}$$

Hence,

$\|\omega(x,s)\| \le \varepsilon^{(1/2)} ( \int\limits_{B(x,s)} \|\omega\|^2 dy)^{(1/2)} <= \varepsilon^{(1/2)} \|\omega(s)\|_{L^2(R^3)}$ Since x is arbitrary, then

$\| \omega(t) \|_{L_\infty(R^3)} \le \varepsilon^{(-1/2)} \| \omega(t) \|_{L^2(R^3)}$ *which is (eq 4.1)* .

Then integrating in time and aplying Cauchy-Schwartz in s we get

$$\int_0^T ||\omega(.,\tau)||_{L^\infty} \, d\tau \le \varepsilon^{(-1/2)} \int\limits_0^t ((\int\limits_{R^3} \| \omega(s,\tau) \|^2 \, ds)^{(1/2)}) d\tau = \varepsilon^{(-1/2)} t^{(1/2)} (\int\limits_0^t (\int\limits_{R^3} \|\omega(x,s) \|^2 \, dx) ds < +\infty$$

Which is the *(eq. 4.2)*

But the right hand side , is the (eq. 4.0), thus it is finite. **QED**

We are now in a position to prove the Conjectures (A) and (B) , non-periodic and periodic setting , homogeneous case of the Millennium problem.

**Corollary 4.1 (Solution of the homogeneous Millennium problem, Conjectures A, B under volume preserving particle structure conservation.)** *For non-compressible flows with the initial data hypotheses as in the formulation of the Millennium problem with finite energy, and without external forces (homogeneous*

*case) , either in the non-periodic or periodic case, and with the physical hypothesis of volume preserving particle structure conservation as in Definition 3.0 ,there cannot occur a blow-up in finite time intervals.*

**Proof:** Immediate after the Theorem 4.1 and the Beale-Kato-Majda theorem 2.2

The Beale-Kato-Majda theorem 2.2 guarantees that if $\int_0^T |\omega(.,\tau)|_{L^\infty}\, d\tau \leq M_1$ in the maximal finite interval [0,T*) that we can have a smooth solution , then we can extend the smooth solution on [0,T*], with no-blowup at T*, thus the unique solution is regular (no- finite time blowup) QED.

## 5. Epilogue

Here is a final , qualitative explanation of the no-blow up solution of the millennium problem will be enlightening, for the reader beyond the formal details.

When the solutions of the millennium problem, in the viscous case satisfy the absolute physical principle of particle structure conservation , there is no possibility, due to non-compressibility, to convert heat to kinetic energy through detonating vortex cones and motions and convert it to kinetic energy. If the fluid was not non-compressible but compressible, this would be possible and the infinite energy reservoir of heat could become an infinite kinetic (rotational) energy situation, thus a blow-up and vortex singularity. But for non-compressible flows that satisfy the physical principle of particle structure conservation this is impossible. We notice that even if for the compressible fluids , it is satisfied a particle structure (which is natural due to the atomic structure of the matter) but not volume preserving particle structure conservation since the volume of the blots at the scale of smalness ε, is not preserved in the definition of the particle structure, we cannot repeat the proof of theorem 4.1, and we cannot prove regularity.

---

# CHAPTER 7

# USE OF THE CKN THEORY TO PROVE THE NON-EXISTENCE OF FINITE TIME VELOCITIES BLOW-UPS TO INFINITE (SINGULARITIES) OF STRONG SOLUTIONS OF NAVIER-STOKES EQUATIONS.

**Konstantine (or Constantine) Kyritsis***

* Associate Prof. of University of Ioannina, Greece. Dept. Accounting-Finance, Psathaki Preveza 48100, ckiritsi@uoi.gr

**Abstract** *In this article, we give a very short proof (Theorem 3.4) that, under the hypotheses of the initial conditions of the millennium problem about the Navier-Stokes equations, there cannot be shaped in finite time, velocities blow-ups to infinite (singularities). We utilize the existing CKN theory of possible blow-ups, in particular, type I and type II, and we prove by reduction to contradiction that neither of the two types of blow-ups can appear in finite time blow-ups because this would mean that a blow-up or singularity will happen in a 3D region too. The proof is based on the non-compressible of the fluid, which implies that the 1D integrals on paths of the pressure force densities is path-independent (conservative vector field) and that the initial energy is finite.*

---

## 1. Introduction

The Clay millennium problem about the Navier-Stokes equations is one of the 7 famous problem of mathematics that the Clay Mathematical Institute has set a high monetary award for its solution. It is considered a difficult problem as it has resisted solving it for almost a whole century. The Navier-Stokes equations are the equations that are considered to govern the flow of fluids, and had been formulated long ago in mathematical physics before it was known that matter consists from atoms. So actually they formulate the old infinite divisible material fluids. Although it is known that under its assumptions of the millennium problem the Navier-Stokes equations have a unique smooth and local in time solution, it was not known if this solution can be extended smoothly and globally for all times, which would be called the **regularity of the Navier-Stokes equations in 3 dimensions.** The corresponding case of regularity in 2 dimensions has long ago been proved to hold but the 3-dimensions had resisted proving it. There are mainly two groups of people with different expectations about this millennium problem. At first it is the group of mathematical physicists, who consider it natural that there should not exist a blow-up in finite time in 3 dimesions either. Then there is the group of mathematicians of partial differential equations, who think that it is very probable that it might exist a blow up in finite time. Both groups seem to be right! The reason is the next. There are physical aspects of the Millennium Problem, that are not formulated mathematically, in the setting of the problem. For example, the energy conservation of the heat of the fluid, is not included in the formulation of the problem at all. There is no symbol or magnitude for the temperature in the formulation of the problem. Still a physicist thinks, that the conservation of the energy of the heat is something that does hold and can use it in the arguments. Nevertheless, a mathematician, considers the millennium problem, as a mathematical problem of partial differential equations, and considers as holding, only what is mathematically stated in the formulation. I belong to the group of mathematical physicists. From the CKN theory of blow-ups in finite time (see [9].[12],[15] we know that blow ups if they happen at all, they will happen on a set of points of Lebesgue measure zero in space and time. *Although the energy and the momentum on a finite 3D region, is always finite, in the relevant integral mathematically we cannot exclude the possibility that the integrated function (energy or momentum) is infinite on a point. The argument in the theorem 3.4, gives a method based on the 1D integral of pressure force densities, to extend this point-infinity (singularity) , on a 1D path and eventually, to a 3D area of blow-up , which will give a contradiction due to finite energy. Thus no point blowup-singularity can exist.*

## 2. The formulation of the millennium problem and the 4 sub-problems (A) , (B), (C), (D)

In this paragraph we highlight the basic parts of the standard formulation of the 4th Clay millennium problem as in [1] Fefferman C.L. 2006.

The **Navier-Stokes** equations are given by (by R we denote the field of the real numbers, ν>0 is the density normalized viscosity coefficient )

$$\frac{\partial}{\partial t}u_i + \sum_{j=1}^{n} u_j \frac{\partial u_i}{\partial x_j} = -\frac{\partial p}{\partial x_i} + \nu \Delta u_i \qquad (x \in R^3 , t>=0 , n=3) \qquad (eq.2.1 )$$

$$div u = \sum_{i=1}^{n} \frac{\partial u_i}{\partial x_i} = 0$$

(x∈$R^3$ , t>=0 , $n$=3) (eq. 2.2 )

with initial conditions $u(x,0)=u^0(x)$ **x∈$R^3$**

and $u^0$ (x) ∈ $C^\infty$ divergence-free vector field on **$R^3$** (eq. 2.3)

If ν=0 then we are taking about the **Euler** equations and **inviscid** case.

$$\Delta = \sum_{i=1}^{n} \frac{\partial^2}{\partial x_i^2}$$

is the Laplacian operator .The Euler equations are (eq2.1), (eq 2.2), (eq 2.3) when ν=0.

It is reminded to the reader, that in the equations of Navier-Stokes, as in (eq. 2.1) the density ρ, is constant, it is custom to normalized to 1 and omit it.

For physically meaningful solutions we want to make sure that $u^0(x)$ does not grow large as |x|→∞. This is set by defining $u^0(x)$ , and f(x,t) and called in this paper **(smooth) Schwartz initial conditions** , to satisfy

$|\partial_x^a u^0(x)| \le C_{a,K}(1+|x|)^{-K}$ on **$R^3$** for any $\alpha$ and K (eq.2.4 )

(Schwartz used such functions to define the space of Schwartz distributions)

**Remark 2.1.** It is important to realize that smooth Schwartz initial velocities after

(eq 2.4) will give that the initial vorticity $\omega_0$ =curl($u^0$) , in its supremum norm, is bounded over all 3-space.

We also assume **external forces** that are Schwartz in space and time:

$|\partial_x^a \partial_t^m f \quad (x,t)| \le C_{a,m,K}(1+|x|+t)^{-K}$on **$R^3$**× $[\mathbf{0},+\infty)$for any $\alpha,m,$K (eq.2.5 )

We accept as physical meaningful solutions only if it satisfies

p, u ∈$C^\infty(R^3\times[0,\infty))$ (eq.2.6 )

and

$\int_{\mathfrak{R}^3} |u(x,t)|^2 \, dx < C$ for all t>=0 (**Bounded or finite energy**) (eq.2.7 )

**Lemma 2.0.** *from the smooth Schwartz initial data of the formulation of the 4th Clay Millennium problem we derive that the supremum of the vorticity at t=0, symbolized by $F_\omega$ exist and it is finite.*

**Proof.** We use the remark 2.1 and of course the inequality |a-b|<=|a|+|b|, at the definition of the vorticity by the partial derivatives. QED.

**Remark 2.2** It is important to realize that smooth external force (densities) with the Schwartz property as in (eq.2.5) , have not only a rule for upper bounded spatial partial derivatives but also the same rule for time upper bounded partial derivatives.

**Remark 2.3** We must stress here that imposing smoothness of the coordinate functions of velocities and external forces of the initial t=0 data and later time t data in the Cartesian coordinates plus and Schwartz condition as in (eq 2.5) is not equivalent with imposing similar such smoothness of the coordinate functions and conditions in the cylindrical or spherical coordinates. Alternatively, to rule out problems at infinity, we may look for spatially periodic solutions of (2.1), (2.2), (2.3). Thus we assume that $u^0(x)$ , and f(x,t) satisfy

$u^0(x+e_j)= u^0(x)$, $f(x+e_j,t)= f(x,t)$, $p(x+e_j,0)=p(x,0)$, for 1<=j<=3 (eq.2.8 )

($e_j$ is the jth unit vector in $R^3$)

In place of (2.4) and (2.5), we assume that $u^0(x)$, is smooth and that

$|\partial_x^a \partial_t^m f \quad (x,t)| \leq C_{a,m,K}(1+t)^{-K}$ on $\mathbf{R^3} \times [\mathbf{0}, +\infty)$ for any *α*,*m*,K (eq.2.9 )

We then accept a solution of (2.1), (2.2) , (2.3) as physically reasonable if it satisfies

$u(x+e_j,t)= u(x, t)$, $p(x+e_j, t)=p(x,t)$, on $\mathbf{R^3} \times [\mathbf{0}, +\infty)$ for 1<=j<=3 (eq.2.10 )

and p, u $\in C^\infty(R^3 \times [0,\infty))$ (eq.2.11 )

In the next paragraphs we may also write $u_0$ instead of $u^0$ for the initial data velocity.

We denote Euclidean balls by $B(a,r) := \{x \in R^3 : ||x-a|| \leq r\}$ , where |x| or ||x|| is the Euclidean norm.

**The 4 sub-problems or conjectures of the millennium problem are the next:**

**(Conjecture A) Existence and smoothness of Navier-Stokes solution on $R^3$.**

*Take v>0 and n=3. Let $u_0(x)$ be any smooth, divergent-free vector field satisfying (2.4). Take f(x,t) to be identically zero. Then there exist smooth functions p(x,t) , u(x,t) on $R^3x[0,+\infty)$ that satisfy (2.1), (2.2), (2.3) , (2.6) , (2.7).*

**(Conjecture B) Existence and smoothness of Navier-Stokes solution on $R^3/Z^3$.**

*Take v>0 and n=3. Let $u_0(x)$ be any smooth, divergent-free vector field satisfying (2.8); we take f(x,t) to be identically zero. Then there exist smooth functions p(x,t) , u(x,t) on $R^3x[0,+\infty)$ that satisfy (2.1), (2.2), (2.3) , (2.10) , (2.11).*

**(Conjecture C) Breakdown of Navier-Stokes solution on $R^3$**

*Take v>0 and n=3. Then there exist a smooth, divergent-free vector field $u_0(x)$ on $R^3$ and a smooth f(x,t) on $R^3x[0,+\infty)$ satisfying (2.4), (2.5) for which there exist no smooth solution (p(x,t) ,u(x,t)) of (2.1), (2.2), (2.3) , (2.6) , (2.7) on $R^3x[0,+\infty)$.*

**(Conjecture D) Breakdown of Navier-Stokes solution on $R^3/Z^3$**

*Take v>0 and n=3. Then there exist a smooth, divergent-free vector field $u_0(x)$ on $R^3$ and a smooth f(x,t) on $R^3x[0,+\infty)$ satisfying (2.8), (2.9) for which there exist no smooth solution (p(x,t) ,u(x,t)) of (2.1), (2.2), (2.3) , (2.10) , (2.11) on $R^3x[0,+\infty)$.*

**Remark 2.4.** It is stated in the same formal formulation of the Clay millennium problem by C. L. Fefferman see [1] Fefferman C.L. 2006 (see page 2nd line 5 from below) that the conjecture (A) has been proved to holds locally. "..if the time internal $[0,\infty)$, is replaced by a small time interval [0,T), with T depending on the initial data....". In other words there is $\infty$>T>0, such that there exists a unique and smooth solution $u(x,t) \in C^{\infty}(R^3 \times [0,T))$. See also [3] A.J. Majda-A.L. Bertozzi ,Theorem 3.4 pp 104. In this paper, as it is standard almost everywhere, the term smooth refers to the space $C^{\infty}$

In the next the $|| \, ||_m$ is the corresponding Sobolev spaces norm and . We denote by $V^m$ ={u in $H^m(R^n)$ and divu=0} where $H^m(R^n)$ are the Sobolev spaces with the $L^2$ norm.

We must mention that in A.J. Majda-A.L. Bertozzi [3] ,Theorem 3.4 pp 104, Local in Time existence of Solutions to the Euler and Navier-Stokes equations it is proved that indeed if the initial velocities belong to $V^m$ m>=[3/2]+2 there exist unique smooth solutions locally in time [0,t]. Here, in the formulation of the millennium problem the hypotheses of smooth with Schwartz condition initial velocities satisfies this condition therefore we have the existence and uniqueness of smooth solution locally in time, both in the non-periodic and the periodic setting without external forcing (homogeneous case).

The existence and uniqueness of a smooth solutions locally in time is stated in the formulation by C.L. Fefferman [1] for the homogeneous cases and conjectures (A), (B). When a smooth Schwartz condition external force is added (inhomogeneous case) , it is natural to expect that also there should exist a local in time unique sooth solution. **And this is indeed provable with the existing bibliography** (The standard theorems are due to **Fujita–Kato**, **Kato**, **Leray**, etc and presented in many textbooks such as Temam [6], Ladyzhenskaya [7] , Majda [3] etc ).

**Lemma 2.1** *The smooth Schwartz initial data and the rest of the hypotheses in the formulation of the $4^{th}$ Cay Millennium problem for the homogeneous case and conjectures (A), (B), and non-homogeneous case and conjectures (C), (D), guarantee, that it exists a unique smooth solution locally in time [0,t) , for the Navier-Stokes equations (without or with external forces).*

We state here also two, very well-known criteria of no blow-up and regularity.

In the next theorem the $|| \, ||_m$ is the corresponding Sobolev spaces norm and . We denote by $V^m$ ={u in $H^m(R^n)$ and divu=0} where $H^m(R^n)$ are the Sobolev spaces with the $L^2$ norm.

**Theorem 2.1 Velocities Sobolev norm sufficient condition of regularity.** *Given an initial condition $u_0 \in V^m$ m>=[3/2]+2=3.5 e.g. m=4 , then for any viscosity v>=0 . there exists a maximal time T* (possibly infinite) of existence of a unique smooth solution $u \in C([0,T^*];V^m) \cap C^1([0,T^*];V^{m-2})$ to the Euler or Navier-Stokes equation. Moreover, if $T^*<+\infty$ then necessarily $lim_{t\to T^*} ||u(.,t)||_m = +\infty$.*

**Proof:** See A.J. Majda-A.L. Bertozzi [3], Corollary 3.2 pp 112). QED

**Remark 2.5** Obviously this proposition covers the periodic case too.

**Theorem 2.2 BKM criterion (Bealo-Kato-Majda) Supremum of vorticity sufficient condition of regularity**

*Let the initial velocity $u_0 \in V^m$ m>=[3/2]+2 , e.g. m=4, so that there exists a classical solution $u \in C^1([0,T] ; C^2 \cap V^m)$ to the 3D Euler or Navier-Stokes equations. Then :*

*(i) If for any T>0 there is $M_1$ >0 such that the vorticity ω=curl(u) satisfies*

$$\int_0^T |\omega(.,\tau)|_{L^\infty} \, d\tau \leq M_1$$

*Then the solution u exists globally in time, $u \in C^1([0,+\infty] ; C^2 \cap V^m)$*

*(ii) If the maximal time $T^*$ or T of the existence of the solution $u \in C^1([0,T] ; C^2 \cap V^m)$ is finite,*

*then necessarily the vorticity accumulates so rapidly that*

$$lim_{t\to T^*} \int_0^T |\omega(.,\tau)|_{L^\infty} \, d\tau = +\infty \qquad (eq. 2.12)$$

**Proof:** See A.J. Majda-A.L. Bertozzi [3], Theorem 3.6 pp 115, $L^\infty$ vorticity control of regularity. QED.

**Remark 2.6** Obviously this proposition covers the periodic case too. But furthermore it is important to know, that the above Bealo-Kato-Majda, holds not only for the **Navier-Stokes** but the for the **Euler** equations and furthermore **without or with external forcing under the Schwartz boundedness conditions above.**

**3. Solution by contradiction on the existence of finite time blow-ups/ singularities (of type I or type II) as in the CKN theory.**

The CKN theory was initially developed by Leray and Caffarelli, L., Kohn, R., & Nirenberg, L. (1982) see [9], [12])But there are other modern expositions of it e.g.

**Lemma 3.1 . (From the CNK theory )** *Any finite time blow-up of velocities, after the standard hypotheses of the formulation of the $4^{th}$ Millennium problem about the Navier-Stokes equations, cannot have singularities on all points of an open 3D , 2D or 1D region. Any possible singularities must have total measure zero.*

**Proof:** (See [9], [12])

**Lemma 3.2. (Galdi G.P.)** *It is impossible that under the standard hypotheses of the formulation of the $4^{th}$ Millennium problem about the Navier-Stokes equations, that there will be created in finite time by blow-up a singularity of the velocities at the spatial infinite.*

**Proof: (See [13] ) .** We may as well formulate the proof for the periodic case , which has compact fundamental region, so as to avoid reasoning, about blow-ups at the infinite. **QED**

**Lemma 3.3.** ***1D INTEGRALS OF PRESSURE FORCES WILL BLOW-UP IN A FINITE TIME VELOCITIES BLOW-UP.*** *Let any point S of velocity blow-up, created in finite time T* after the standard hypotheses of the formulation of the $4^{th}$ Millennium problem about the Navier-Stokes equations. Then at time T*-δt*

*For a very small time δt , (almost infinitesimal) just before the blow up, any path P=AS from a fixed in time and space point A, till the point S will define an 1D integral of the work density of the pressure forces densities, which becomes infinite as δt converges to 0, and the same even if we weight the integral in polar coordinates by $r^2$, where r=||x-S||.*

**Proof:** Direct from CKN theory of possible finite time ,blow-ups (see [12]).

In the classical work by Caffarelli R, Kohn R, Nirenberg L, "Partial Regularity of suitable Weak Solutions of the Navier-Stokes Equations" See [12] that is known by now as the CKN theory, in page 795, paragraph 5, under the title "**The Blow-up estimate**" they write that ".....is to view Corollary 1, as an estimate for **the** **minimum rate** at which a singularity can develop....that near a singular point ($x_0$,$t_0$) $||u||>= C/r$ where $r^2=||x-x_0||^2 +||t-t_0||^2$ ......." where u is the velocity. And ||,|| is the euclidean norm.

Later authors classified as singularities-blow-ups of **type I**, those that have this minimum speed, and of **type II**, if they have higher speed of divergence of the velocities.

In our symbolism here $S =x_0$ is a point of blow=up at time $T^*=t_0$, r=R

In the CKN theory of singularities **s**ee e.g. [15] In page 9 Proposition 2, the author remarks : (1) "*This means that the only way a singularity can manifest itself, in a Leray solution of the NS equations, is through a divergence (blowup) of the velocity field itself*." A strong solutions after the initial conditions of the millennium problem, as is our case, is a special case of a Leray-Hopf weak solution. Thus whenever a blow-up (singularity) occurs in finite time, in the sense of breaking the smoothness (thus a vorticity accumulation blow-up after the BKM criterion in theorem 2.2) , the blow up must be also manifest at the velocities too.

In figuring out if given the hypotheses of the initial conditions of the millennium problem, a blowup or not (regularity) in finite time of a physical quantity (velocities, pressures, vorticity, accelerations, pressure forces etc.) will imply a blowup or not (regularity) in finite time of another quantity a general rule is the next:

*a) Ascending the order of derivatives, propagates the implication of blowups*

*b) Descending the orders of derivatives (integration) propagates the implication of regularity.*

Following this rule and ordinary differential calculus , it is easy to verify that once the velocity blows-up , also the pressures blow-up and also the pressure forces blow-up too.

We will consider that our estimates of speed of divergence, are a very small time δt (almost infinitesimal) just before the blow-up time T* which is constant though out the **spatial limits and allows us to deal with strong solutions**.

Thus for a sufficient small ball around S, and almost every where

$||u||>= C/R$ where $R^2=||x-S||^2 +||\delta t||^2$ , $r=||x-S||$

We could as well use the omega symbolism

$||u||=\Omega(C/R)$ sufficient close to S, and a.e (eq. 3.4.0)

The relation of pressure and velocities is the well know, by the Riesz transform of the velocity terms

-Δp=$\sum_{ij}\partial_i\partial_j(u_iu_j)$

This implies that p scales as $u^2$, and that pressures will blow up if the velocities will blow up.

It is well known (e.g. in the work of Seregin and Sverak, but also by direct spacial integration) that given the above minimal speed of spatial divergence of the velocity, then the pressure diverges with **spatial speed**

$||p(x, T^*-\delta\tau)|| >= C_2/r^2$, sufficient close to S, and a.e. (eq. 3.4.1)

Then the pressure force densities $F_p$, $||F_p|| = \|-\nabla p\|$ will blow up in spatial divergence by a rate of

$||F_p(x, T^*-\delta\tau)|| >= C_3/r^3$, sufficient close in space time to the singularity to and a.e. (eq. 3.4.2)

Then a 1D path integral $I(F_p)$ of the pressure forces densities sufficient close to T* , will blow up as spatial divergence also at a rate lower by 1, due to integration

$||I(F_p)|| >= C_4/r^2$, sufficient close to the time of singularity , and a.e.

(eq. 3.4.3)

In addition, if e.g. we set polar coordinates (r,θ,ϕ) , with center, the singularity point S and a ball of radius R around it , then an integral of the **pressure forces densities** weighted by $r^2$ , the corresponding integral of **work densities** of the force densities along a radius, will also blow-up , sufficient close in space time to the singularity , and a.e.

$|I_2(F_p)| = \| \int_S^R (F_p(r,\theta,\phi), n) r^2 dr \|$ >=C5|Log(r)| , sufficient close to the singularity time, and a.e.

(eq. 3.4.4)

As δt converges to 0, the above integral blows-up. The $(F_p, n)$ is the inner product of the pressure **force density**, with unit normal vector along the radius path AS. Then the integral is interpreted as a **work-density along the a radius AS.**

The reason is that the $F_p$ diverges with speed $C_3/r^3$ and after the weight $r^2$ , it will diverge with speed $C_6/r$ inside the integral, thus the integral will diverge as logarithm at zero. See also [9] [10] ,[12].[13], [14]. QED.

**THEOREM 3.4. THE CRUCIAL THEOREM FOR THE SOLUTION.** *There cannot be created any finite point S of blow-up of velocities, within a finite time interval, after the standard hypotheses of the formulation of the 4th Millennium problem about the Navier-Stokes equations, as this would lead that it must be necessarily a singularity and blow-up on an open 3D region of S too!*

**Proof:** When searching for a proof that under the initial conditions of the millennium problem, there cannot be shaped a finite time velocity blow-up, the idea came to me that if I could prove that the a 3D volume integral of the pressure force densities , that represents total work of them, around the singularity will blow up, then this would lead to a contradiction, because of the finite energy at the initial conditions, and also because the CKN theory (again by reasons of energy) cannot allow a whole 3D region to blow up.

The proof will by reduction to contradiction. Let us assume that at finite time T*, with the standard hypotheses of the Millennium problem , about the Navier-Stokes equations, a velocity blow-up singularity S is shaped (of any type , either I or I etc) Let a ball B(S, R) with center the singularity S, and of radius R sufficient small, so that the minimal space-time speed of blow ups is definable and exist as in the Lemma 3.3 above . And let Ar(S,R)= $\partial B(S,R)$, be the 2D spherical surface-boundary of the ball.

We consider polar coordinates (r,θ,ϕ) with center the singularity S as paramemetrization of integrals.

Now integrate the pressure force densities , with this parametrization (radii first in the inner integral then ball surface in the outer integral ) the pressure **forces densities** $F_p$, we will get an integral in polar coordinates interpreted as integral of **work densities**, of the form:

$$\iiint_B (F_p, \hat{n})dV = \iint_{A\in\partial B} \int_A^S (F_p(r,\theta,\phi), \hat{n})r^2 \sin(\theta)dr)d\theta d\phi \qquad \text{(eq. 3.4.5)}$$

Where $dV = r^2 sin(\theta)drd\theta d\varphi$

Because we know from the proof of the previous lemma 3.3, and in equation (eq. 3.4.4) , that the most inside 1D path integral in the integral iteration, will blow up as δ->0 for all points A, then we conclude that this will happen for the 3D integral as well.

$$\iint_{A\in\partial B} (\int_A^S (F_p(r,\theta,\phi), \hat{n})r^2 \sin(\theta)dr)d\theta d\phi = \iiint_B (F_p, \hat{n})dV = +\infty \qquad \text{(eq. 3.4.6)}$$

Where $dV = r^2 sin(\theta) dr d\theta d\varphi$

Or more exactly

$$\int_0^{2\pi}\int_0^{\pi}\int_0^{R} (F_p(r,\theta,\phi),\hat{n}) r^2 \sin(\theta) dr d\theta d\phi = +\infty \qquad \text{(eq. 3.4.7)}$$

Nevertheless this integral **has the dimensions of energy, or work**, because the $F_p$ is a 3D force density, the most inside integral a work density of 3D force density, and thus the overall 3D integral will have the dimensions of work or energy. The interpretation is of an "instanteneous" **work by the pressure forces**, or as some people would say the **work required to turn the pressures from the 2D surface of the ball B, to the pressures of the point S** (or if at δt before the blow-up, to the pressures of the surface of a small ball aroundthe point S at time T*-δt).

**This work or energy being infinite, is impossible, as it we have only finite initial energy!**

**We conclude therefore, by reduction to contradiction, that there cannot exist a point-singularity** (of any type, type I or type II etc). QED.

**4. The Solution of the 4th Millennium problem about the Navier-Stokes, in favor of regularity (no singularities, Conjecture A, B).**

Since it has been proved by other authors of CKN theory of blow-ups that, the set of blow-up singularities in a finite time blow-up ,

a) must have measure zero (in space and time),

b) it does not include the infinite point, (or we formulate them for the periodic case)

c) and that all velocity blow-up singularities are of type I or type II,

then, the previous lemma 3.4 that forbids also velocity blow-ups of type I, and type II, proves that there cannot exist any velocity blow-ups at all by finite time blow-up. Which seems to prove the Conjecture A or B which is would become a theorem.

**Nevertheless at this point, we must mention a subtle point.** The formulation of the 4th Millennium problem by C.L. Fefferman as in [1], is stated such that the negation of the **BKM criterion** of regularity (as in Theorem 2.2 ) is that there is a **velocity blow-up** , that as the term suggest, the norms of some velocities **converge to +∞**, and therefore the smoothness and continuity (regularity) also stops. This is the understanding of the Fefferman in his presentation of the formulation of the problem. We will call it **Fefferman implicit hypothesis of velocities blow-up (in other words that a blowup in finite time of the vorticity accumulation, will imply a blowup of the velocities in finite time).** Nevertheless strictly logically, the negation of the **BKM criterion** of regularity (as in Theorem 2.2 ) is that there is a **velocity regularity stop or break** (in finite time). I am not aware of any published proof,**besides the CKN theory** , **that a finite time velocity blow-up** is equivalent to, **a finite time velocity regularity stop!** Obviously of course the first implies the second. In figuring out if given the hypotheses of the initial conditions of the millennium problem, a blowup or not (regularity) in finite time of a physical quantity (velocities, pressures, vorticity, accelerations, pressure forces etc.) will imply a blowup or not (regularity) in finite time of another quantity a general rule is the next:

*a) Ascending the order of derivatives, propagates the implication of blowups*

*b) Descending the orders of derivatives (integration) propagates the implication of regularity.*

Thus a blowup in finite time of the integral of vorticity accumulation (see Theorem 2.2) will imply a blow-up of the vorticity. But since ω=curl(u), and the operator curl involves $1^{st}$ order of spatial derivatives of the velocity u , we cannot directly imply with this rule that a finite time blow-up of the vorticity ω, implies a finite time blowup of the velocities.

Still, the fact that we are in mathematical physics rather, than pure mathematics, gives a way out. The physical magnitude dimensions of the vorticity is that of **circular frequency** ω=2πν, where ν is the frequency, and π is the irrational number 3.14….And in circular motion it is hard to imagine that when the frequency of rotation, goes to infinite (blows-up) , the velocities also do not blowup. Actually the relation of the circular velocity ω, with the linear velocity, **in circular motion** is the ω=u/r . So for a very small **but fixed finite radius Δr** from the center of rotation , from ωΔr=u, **we deduce that when ω blows up, the velocity also must blow up.** This suggest, that **from the physical point of view, the Fefferman implicit hypothesis is valid.**

A full proof of the Fefferman implicit hypothesis comes from the CKN theory of singularities See e.g. [15] In page 9 Proposition 2, the author remarks : (1) "*This means that the only way a singularity can manifest itself, in a Leray solution of the NS equations, is through a divergence (blowup) of the velocity field itself.*" A strong solutions after the initial conditions of the millennium problem, as is our case, is a special case of a Leray-Hopf weak solution, thus we have that the Fefferman implicit hypothesis has been proved in the CKN theory of singularities.

Strictly speaking **we have just proved that under the hypotheses of the initial conditions of the millennium problem, and if the CKN theory is valid, that there cannot be shaped in finite time any velocity blow-ups.** Of course under the Fefferman implicit hypothesis which has been proved in CKN theory, this also means that the velocities regularity cannot stop in finite time, **which solves the millennium problem.**

**(THEOREM A) Existence and smoothness of Navier-Stokes solution on $R^3$ ()**

*Take v>0 and n=3. Let $u_0(x)$ be any smooth, divergent-free vector field satisfying (2.4). Take f(x,t) to be identically zero. Then there exist smooth functions p(x,t) , u(x,t) on $R^3x[0,+\infty)$ that satisfy (2.1), (2.2), (2.3) , (2.6) , (2.7).*

**(THEOREM B) Existence and smoothness of Navier-Stokes solution on $R^3/Z^3$.**

*Take v>0 and n=3. Let $u_0(x)$ be any smooth, divergent-free vector field satisfying (2.8); we take f(x,t) to be identically zero. Then there exist smooth functions p(x,t) , u(x,t) on $R^3x[0,+\infty)$ that satisfy (2.1), (2.2), (2.3) , (2.10) , (2.11).*

**5. Epilogue.** Compressible fluids under the standard hypotheses of the $4^{th}$ millennium problem about the Navier-Stokes equations , do give cases of blow-up in finite time. The reason is that tornadoes can be shaped that stretch the vorticity, and by detonation extract infinite energy from the heat energy of the (infinite space) fluid, thus creating a finite time blow-up and singularities. But non-compressible fluids with finite initial energy cannot do that, thus there are not finite time blow-ups, and the flow is smooth and regular till infinite time.

---